\documentclass[a4paper,11pt]{article}
\usepackage{amssymb}
\usepackage{cite}

\newcommand{\s}{s}

\newcommand{\cL}{{\cal L}}
\newcommand{\cK}{{\cal K}}

\newcommand{\R}{\mathbb{R}}

\newcommand{\N}{\mathbb{N}}

\newcommand{\cF}{{\cal F}}
\newcommand{\cO}{{\cal O}}
\newcommand{\cD}{{\cal D}}
\newcommand{\cuad}{{\sqcap\kern-.68em\sqcup}}

\newtheorem{definition}{Definition}[section]
\newtheorem{theorem}{Theorem}[section]
\newtheorem{proposition}{Proposition}[section]

\newtheorem{lemma}{Lemma}[section]
\newtheorem{corollary}{Corollary}[section]
\newtheorem{remark}{Remark}[section]
\newcommand{\bremark}{\begin{remark} \em}
\newcommand{\eremark}{\end{remark} }

\begin{document}

\begin{center}{\bf      Liouville Theorem for semilinear  elliptic  inequalities     \\[1.5mm]

with the fractional  Hardy operators

 }\bigskip \bigskip

 {\small    Huyuan Chen\footnote{chenhuyuan@yeah.net} \qquad Ying Wang\footnote{ yingwang00@126.com}
\medskip

   School of Mathematics and Statistics, Jiangxi Normal University, Nanchang,\\
Jiangxi 330022, PR China \\[4mm]

 Hichem Hajaiej\footnote{  hichem.hajaiej@gmail.com}
\medskip

California State University, Los Angeles, 5151, USA\\[8mm]

}

\begin{abstract}
In this paper,  we  give a full classification of the nonexistence of positive weak solutions to the semilinear   elliptic  inequality  involving the  fractional Hardy potential
$$(E) \qquad\qquad\qquad\qquad\quad  (-\Delta)^s u+\frac{\mu}{|x|^{2s}} u\geq Q u^p\qquad\qquad\qquad \qquad   $$
in a bounded punctured domain   or in an exterior domain, where $p>0$, $\mu\geq\mu_0$ and
$\liminf Q(x)|x|^{-\theta}>0$ for some $\theta\in\R$,   $-\mu_0$ is the best constant in the fractional Hardy inequality.  Our work covers the important case $\theta\leq -2s$ for which all the existence methods do not apply to get the Liouville type theorems.

In the punctured domain, Br\'{e}zis-Dupaigne-Tesei proved the nonexistence of positive solutions of $(E)$ when $s=1$, $\mu\in[\mu_0,0)$ and $Q=1$. When $s\in(0,1)$, we extend this nonexistence result to a general setting via a new method. We also provide a critical exponent $p^\#_{\mu, \theta}$ depending on the parameters $\mu\in[\mu_0,0)$, and $\theta\in\R$. In the super-critical case, $p>p^\#_{\mu, \theta},$ we obtain the nonexistence of positive solutions for $(E)$,     the nonexistence in the critical case $p=p^\#_{\mu, \theta}$ holds true if $p^\#_{\mu, \theta}>1$.  When $p=p^\#_{\mu, \theta}\in(0,1),$ positive solutions can be constructed  provided that we have an appropriate upper bound of $Q$ near the origin additionally.

In an exterior domain,  we provide the Serrin's type critical exponent $p^*_{\mu, \theta}$ depending on the parameters $\mu\geq \mu_0$,  $\theta\in\R$, and then we obtain the nonexistence of positive solutions in the subcritical case: $p\in(0,p^*_{\mu, \theta})$, the nonexistence  in the critical case $p=p^*_{\mu, \theta}>1$ and the existence of positive solution in the critical case  $p=p^*_{\mu, \theta}\in(0,1)$ under additional assumptions on the upper bound of $Q$ at infinity.

Our methods are self-contained and new. The main ideas and key ingredients will be discussed in the next section after Theorem \ref{teo 1} for punctured domains, and after Theorem \ref{teo 2} for exterior domains. We will also explain why all the previous methods and techniques do not apply to our general setting. Let us give here a foretaste: Based on  the imbalance between  the Hardy operator and the nonlinearity,   we can obtain an initial asymptotic behavior rate at the origin or at infinity, we then  improve this rate  by using the interaction with the nonlinearity. By repeating this process finite number of times,  a contradiction will be deduced from  the nonexistence for the related non-homogeneous fractional Hardy problem. This process allows us to obtain    the  nonexistence for the fractional Hardy problem with  larger ranges of $\theta$  and $p$. Our study covers all possible ranges, and our results are optimal.

\end{abstract}

  \end{center}
% \tableofcontents \vspace{1mm}
  \noindent {\small {\bf Keywords}:  Liouville Theorem; Fractional Hardy operator.  }
   \smallskip

   \noindent {\small {\bf AMS Subject Classifications}:     35R11;  35J75; 35A01.
  \smallskip

\vspace{2mm}

\setcounter{equation}{0}
\section{Introduction}
 Let $\Omega$ be a bounded   $C^2$ domain containing the origin  in $\R^N$ with $N\geq 2$, and the fractional Hardy operator be defined by
$$ \cL_\mu^s  = (-\Delta)^s    +\frac{\mu}{|x|^{2s}} $$
for $\s\in(0,1)$   and $$\mu\geq  \mu_0 :=-2^{2\s}\frac{\Gamma^2(\frac{N+2\s}4)}{\Gamma^2(\frac{N-2\s}{4})},$$
where  $ (-\Delta)^s $ is the fractional Laplacian  defined by
$$
(-\Delta)^\s  u(x)=\lim_{\varepsilon\to0^+} (-\Delta)^\s_\varepsilon   u(x)
$$
here $$ (-\Delta)^\s_\varepsilon   u(x) =C_{N,\s}\int_{\R^N\setminus B_\varepsilon(x) }\frac{ u(x)-
u(z)}{|x-z|^{N+2\s}}  dz,$$  $B_\varepsilon(x) \subset \R^N$ is the ball of radius $\varepsilon>0$ centered at $x$, $B_\varepsilon=B_\varepsilon(0)$,  $C_{N,\s}=2^{2\s}\pi^{-\frac N2}\s\frac{\Gamma(\frac{N+2\s}2)}{\Gamma(1-\s)}$ with $\Gamma$ being the Gamma function.

The main goal of this paper is to discuss  the nonexistence of positive weak solutions to the elliptic inequality
\begin{equation}\label{eq 1.1}
 \cL_\mu^s u\geq Q   u^p \quad {\rm in}\ \,  \cO,
 \end{equation}
subject to the Dirichlet boundary type condition: $u\geq0$ in $\R^N\setminus \Omega$ when $\cO=\Omega\setminus\{0\}$ is a bounded punctured domain
 or $u\geq0$ in $\bar \Omega$ in $\cO=\R^N\setminus \bar\Omega$ is an exterior domain,
where   $p>0$
  and the potential $Q\in C^\beta_{loc}(\R^N\setminus\{0\})$ with $\beta\in(0,1)$ is a positive function such that for some $\theta\in\R$,
\begin{equation}\label{pt r1}
  \liminf_{|x|\to0^+}Q(x)|x|^{-\theta}>0\quad{\rm  for\ the\ punctured\ domain}
\end{equation}
and
\begin{equation}\label{pt r2}
   \liminf_{|x|\to+\infty}Q(x)|x|^{-\theta}>0\quad{\rm for\ the\ exterior\ domain}.
\end{equation}
  Let
$$L^1_s(\R^N)=\big\{w\in L^1_{loc}(\R^N):\, \int_{\R^N}\frac{|w(x)|}{1+|x|^{N+2\s}}dx<+\infty\big\}.$$
\begin{definition}\label{weak definition}
$~$
We say that $u$ is a classical solution of (\ref{eq 1.1}), if $u\in
L^1_s(\R^N)\cap C(\cO)$ satisfies (\ref{eq 1.1}) pointwisely. \smallskip

$-$ We say that $u$ is a positive solution of (\ref{eq 1.1}), if $u$ is a classical solution of (\ref{eq 1.1}) and $u>0$  a.e. in $\cO$.
 \smallskip

$-$ We say that $u$ is a weak   solution of (\ref{eq 1.1}), if $u\in
L^1_s(\R^N)$, $Q|u|^p\in L^1_{loc}(\cO\setminus\{0\})$,  for $\epsilon>0$ small enough,  letting
\begin{equation}\label{kkk}
\cK_\epsilon=\cO\setminus \overline{B_\epsilon}\ \ {\rm if}\ \, \cO=\Omega\qquad{\rm or}\qquad
\cK_\epsilon=\cO\setminus \overline{B_{\frac1\epsilon}}\ \ {\rm if}\ \,\cO=\R^N\setminus \bar \Omega,
\end{equation}
\begin{equation}\label{weak sense}
\int_{\cK_\epsilon}  u \cL_\mu^s\xi \, dx\geq \int_{\cK_\epsilon} Qu^p \xi
dx,\quad \forall\,\xi\in \mathbb{X}_{s,+}(\cK_\epsilon),
\end{equation}
where $\mathbb{X}_{s,+}(\cK_\epsilon)=\big\{\xi\in\mathbb{X}_{s}(\cK_\epsilon), \xi\geq0 \ {\rm in\  } \R^N\big\},$ $\mathbb{X}_{s}(\cK_\epsilon)\subset C(\R^N)$ is the space of   functions
$\xi$ satisfying:  \smallskip

\noindent (i)     $\rm{supp}(\xi)\subset \overline{\cK_\epsilon}$;\smallskip

\noindent(ii) $(-\Delta)^s\xi(x)$ exists for all $x\in \cK_\epsilon$
and $|(-\Delta)^s\xi(x)|\leq C$ for some $C>0$;\smallskip

\noindent(iii) there exist $\varphi\in L^1(\cK,\rho_{_{\cK_\epsilon}}^s dx)$
and $\varepsilon_0>0$ such that $|(-\Delta)_\varepsilon^s\xi|\le
\varphi$ a.e. in $\cK$, for all
$\varepsilon\in(0,\varepsilon_0]$, where
$\rho_{_{\cK_\epsilon}}={\rm dist}(x,\partial \cK_\epsilon)$.
\end{definition}
\begin{remark}
$(i)$ A sub-weak   solution of
 $$\cL_\mu^s u\leq Q   |u|^{p-1}u \quad {\rm in}\ \,  \cO, \qquad u\leq 0\ \, {\rm in}\ \, \R^N\setminus \cO$$
 can be defined as follows:
 If $u\in L^1_s(\R^N)$ is nonpositive, $Q|u|^p\in L^1_{loc}(\cO\setminus\{0\})$ and  for $\epsilon>0$
$$%\begin{equation}\label{weak sense-}
\int_{\cK_\epsilon}  u \cL_\mu^s\xi \, dx\leq \int_{\cK_\epsilon} Q|u|^{p-1}u  \xi
dx,\quad \forall\,\xi\in \mathbb{X}_{s,+}(\cK_\epsilon).
$$%\end{equation}

$(ii)$   Since   the fractional Hardy operator degenerates at the poles at the origin and at infinity,   our definition of weak solution  has to be set to work on bounded domains away from the poles.
 The advantage of our definition of weak solutions is that the Kato's inequality holds true in this setting.

 \end{remark}
  Note that this  definition of weak solutions was introduced in \cite{CV0}, where the authors dealt with
the weak solutions of semilinear fractional problems with Radon measure, subject to
zero Dirichlet boundary condition.

The prototype of (\ref{eq 1.1}) is the Lane-Emden equation
 $$-\Delta u=u^p\quad{\rm in}\ \ \R^N\setminus\{0\},$$
which doesn't admit any  positive solution for $p\in(1, \frac{N}{N-2}]$, see \cite{AS,BFV,CPZ,CL,PT}. The nonexistence in exterior domains was established in \cite{CF}.
When  $p>\frac{N}{N-2},$  it has a family of fast decaying solutions at infinity and one slow decaying solution  by the phase plane analysis \cite{MP,P1}.  Other closely related problems were discussed in \cite{L2,LLd}  and the references therein.  The Hardy operator with $s=1$ has the form:
 $\cL_\mu:=-\Delta+\frac{\mu}{|x|^2}$.  For $N\ge3$ and  $\mu>\frac{(N-2)^2}{4}$, the semilinear Hardy problem
\begin{equation}\label{eq Hardy}
\cL_\mu u=|x|^{\theta}u^p \quad\ {\rm in}\ \  \R^N\setminus\{0\}
\end{equation}
has been studied extensively,  one-parameter family of finite-energy radial solutions was obtained in \cite{BV,T}, and the global regularity was established in \cite{W1}.  In a bounded domain containing the origin, \cite{CZ} studied the isolated singularities of the problem
$$%\begin{equation}\label{eq Hardy punc}
\cL_\mu u=u^p\quad  {\rm in}\ \, \Omega\setminus\{0\}, \qquad u=0\quad{\rm on}\ \, \partial\Omega
$$%\end{equation}
by building suitable connections  with  the weak solutions of
$$
    \displaystyle    \mathcal{L}_\mu    u-u^p =k\delta_0\quad
 {\rm in}\ \, \Omega, \qquad u=0\quad{\rm on}\ \, \partial\Omega
$$
in a  weighted  distributional sense inspired by \cite{CQZ} for
$ p\in(1,1+\frac{2}{\frac{N-2}{2}+\sqrt{ (N-2)^2/4+\mu}})$, where the upper bound is
the corresponding Serrin critical exponent. More related  discussions can be found in \cite{CV1,CV2}.  When  $\theta\in(-2,2)$  and
 $p=1+\frac{2+\theta}{\frac{N-2}{2}+\sqrt{ (N-2)^2/4+\mu}}$. \cite{FG} classified the solutions and established that any solution of  problem (\ref{eq Hardy})  has  the following isolated singularity:
$$\lim_{|x|\to0^+} u_p(x) |x|^{\frac{2+\theta}{p-1}}(-\ln|x|)^{\frac1{p-1}} =c_{p,\mu,\theta},$$
 where $c_{p,\mu,\theta}>0.%=\Big[\frac{2+\theta}{p-1}(N-2-\frac{2+\theta}{p-1})+\mu\Big]^{\frac1{p-1}}.
 $
The authors in \cite{FG,MV,CY} obtained the fast decaying  solution and slow decaying solution
in the supercritical case as well as the uniqueness.  \cite{CF} classified the isolated singular solutions of the Hardy problem involving the absorption nonlinearity.  \smallskip

Motivated by various applications in different fields and developments in the theory of PDEs,  there has been an increasing interest in the Dirichlet problems with nonlocal operators, especially, the fractional Laplacian, see
e.g. \cite{DPV,Ls07,JLX} and the references therein.
For the basic properties of the fractional Laplacian, we refer e.g. to \cite{DPV,RS,RS0}. It is known that  $(-\Delta)^\s  u(x)$ is well-defined if $u$ is twice continuously differentiable in a neighborhood of $x$ and contained in the space $L^1_s(\R^N)$.
    We also note that, for $u \in L^1_s(\R^N)$, the fractional Laplacian $(-\Delta)^\s  u$ can also be defined as a distribution in the following sense:
$$%\begin{equation}
 % \label{eq:-distributional-sense}
 \langle(-\Delta)^\s u, \varphi \rangle= \int_{\R^N} u (-\Delta)^\s \varphi \,dx \quad\  {\rm for\ any }\ \  \varphi \in C^\infty_c(\Omega).
$$%\end{equation}
We then have
$$%\begin{equation}
 % \label{eq:-distributional-sense}
\cF((-\Delta)^\s u) = |\cdot|^{2\s} \hat u \qquad  {\rm in}\ \  \R^N \setminus \{0\}
$$%\end{equation}
in the sense of distributions,  both $\cF(\cdot)$ and $\hat \cdot$ denote the Fourier transform.

Recently,  it was shown in \cite{CW} that  for $\mu\geq \mu_0$ the equation
$$\mathcal{ L}_\mu^s u=0\quad{\rm in}\ \ \R^N\setminus \{0\}$$
 has two distinct radial solutions
 $$\Phi_{s,\mu}(x)=\left\{\arraycolsep=1pt
\begin{array}{lll}
 |x|^{\tau_-(s,\mu)}\quad
   &{\rm if}\quad \mu>\mu_0\\[1.5mm]
 \phantom{   }
|x|^{-\frac{N-2s}{2}}\ln\left(\frac{1}{|x|}\right) \quad  &{\rm   if}\quad \mu=\mu_0
 \end{array}
 \right.\qquad \ {\rm and}\quad\ \ \Gamma_{s,\mu}(x)=|x|^{\tau_+(s,\mu)},$$
where
$\tau_-(s,\mu)  \leq  \tau_+(s,\mu)$ are such that:
 \begin{eqnarray*}
 &\tau_-(\s,\mu)+\tau_+(\s,\mu) =2\s-N \quad\  {\rm for\ all}\ \   \mu \ge \mu_0,\\[1.5mm]
&\tau_-(\s,\mu_0)=\tau_+(\s,\mu_0)=\frac{2\s-N}2,\quad\
\tau_-(\s,0)=2\s-N, \quad\ \tau_+(\s,0)=0,  \\[1.5mm]
&\displaystyle\lim_{\mu\to+\infty} \tau_-(\s,\mu)=-N\quad {\rm and}\quad \lim_{\mu\to+\infty} \tau_+(\s,\mu)=2\s.
 \end{eqnarray*}
In the sequel and when there is no ambiguity,  we denote $\tau_+=\tau_+(s,\mu)$, $\tau_-=\tau_-(s,\mu)$.    For any $\xi\in C^2_c(\R^N)$, we have
$$%\begin{equation}\label{1.2}
 \int_{\R^N}\Phi_{s,\mu}   (-\Delta)^\s_{\Gamma_{s,\mu}}\xi \, dx  =c_{\s,\mu}\xi(0),
$$% \end{equation}
where $c_{\s,\mu}>0$ and  the dual of the operator of $\cL^s_\mu$ is $(-\Delta)^\s_{\Gamma_{s,\mu}}$, which is  a weighted fractional Laplacian  given by
$$%\begin{equation}\label{L}
(-\Delta)^\s_{\Gamma_{s,\mu}} v(x):=
C_{N,\s}\lim_{\epsilon\to0^+} \int_{\R^N\setminus B_\epsilon(x) }\frac{v(x)-
v(z)}{|x-z|^{N+2\s}} \, \Gamma_{s,\mu}(z) dz.
$$%\end{equation}
Via the above weighted distributional form, the isolated singular solutions of the nonhomogeneous Hardy problem
$$
\cL_\mu^s u= f \quad
    {\rm in}\ \ \Omega\setminus\{0\}, \qquad
u\geq0 \quad   {\rm in}\ \  \R^N\setminus  \Omega
$$
have been classified  under an optimal assumption for nonnegative functions $f\in C^\beta_{loc}(\bar \Omega\setminus\{0\})$ with $\beta\in(0,1)$, see Theorem \ref{theorem-C} in Section 2.

Inspired by the above classification of nonhomogeneous problems.   The goal of this paper is to show the nonexistence of positive solutions of fractional Hardy problem (\ref{eq 1.1}) with a source nonlinearity.  Before stating our main results,   let us introduce two  critical exponents
\begin{equation}\label{eq 1.1-cr1}
p^*_{\mu,\theta}=  1+\frac{2s+\theta}{-\tau_-}
\quad\ \ {\rm and}\quad\ \
p^\#_{\mu,\theta}=\left\{
\begin{array}{lll}
  1+\frac{2s+\theta}{-\tau_+} &\quad {\rm if }\;\ \mu_0\leq \mu<0,\\[2mm]
+\infty&\quad  {\rm if }\;\ \mu\geq 0.
\end{array}\right.
\end{equation}
Here $p^*_{\mu,\theta}$ is the related Serrin critical exponent and $p^\#_{\theta,\mu}$
is a particular exponent defined for $\mu\in[\mu_0,0)$.
Moreover, for $\mu\in[\mu_0,0)$,
$$\mathbf{
 p^\#_{\mu,\theta}\geq 1\ \  {\rm if}\ \,  \theta\geq -2s,\quad\  p^\#_{\mu,\theta}<1\ \   {\rm if}\ \,   \theta<-2s,\quad\  p^\#_{\mu,\theta}\leq0 \ \   {\rm if}\ \,   \theta\leq  \tau_+-2s
 }$$
and
$$\mathbf{
p^*_{\mu_0,\theta}=p^\#_{\mu_0,\theta}=\frac{N+2s+2\theta}{N-2s}.
}$$

%When $\mu>0$,  $\tau_+>0$ and
%$$p^\#_{\theta,\mu}\geq 1\ \ \ {\rm if}\ \ \, \theta\leq  -2s.$$

Now we state our  nonexistence result for (\ref{eq 1.1}) in a bounded punctured domain, that is,
\begin{equation}\label{eq 1.1-in}
 \cL_\mu^s u\geq Q    u^p \quad
    {\rm in}\ \  \Omega\setminus\{0\}, \qquad\quad
u\geq0 \quad  {\rm in}\ \ \R^N\setminus \Omega.
\end{equation}

\begin{theorem}\label{teo 1}
Let  $\mu\in[\mu_0, 0)$ and $Q$ verify (\ref{pt r1})  with $\theta\in\R$ and $p^\#_{\mu,\theta}$ be defined in (\ref{eq 1.1-cr1}).\smallskip

$(i)$ If $\mu=\mu_0$,  then problem (\ref{eq 1.1-in})    has no positive weak solution for any $\theta\in\R$ and $p\geq \max\{p^\#_{\mu,\theta},0\}$;

$(ii)$ If $\mu\in(\mu_0, 0)$, $\theta>-2s$ and  $p\geq p^\#_{\mu,\theta}$,  then problem (\ref{eq 1.1-in})    has no positive weak solution;\smallskip

$(iii)$ If $\mu\in(\mu_0, 0)$, $\theta= -2s$, either $p>1$ or
 $$p=1\ \ \  {\rm and}\ \ \  \liminf_{|x|\to0}Q(x)|x|^{2s}>\mu-\mu_0,$$
  then problem (\ref{eq 1.1-in})   has no positive weak solution;\smallskip

  $(iv)$  If $\mu\in(\mu_0, 0)$,  $\theta< -2s$ and  $p> \max\{p^\#_{\mu,\theta},0\}$, then problem (\ref{eq 1.1-in})    has no positive weak solution.

\end{theorem}

We should mention that   the nonexistence of positive solution for the semilinear Hardy inequality
$$%\begin{equation}\label{eq Hardy punc}
\cL_\mu u\geq u^p\quad  {\rm in}\ \, B_r\setminus\{0\}, \qquad u=0\quad{\rm on}\ \, \partial B_r
$$%\end{equation}
 was first studied in \cite{BDT}, more related results can be found in \cite{D-1,FM}. When $s\in (0,1)$ and $\mu\in(\mu_0, 0)$, \cite{F0} established the nonexistence
of  (\ref{eq 1.1-in})  for $ p\geq p^\#_{\mu,0}$ when $Q\equiv1$, via a harmonic extension given by Caffarelli-Silvestre \cite{CS1}. \smallskip

Our derivation of the  nonexistence in a bounded punctured domain comes from the observation that $\Gamma_{s,\mu}$ blows up at the origin when $\mu\in[\mu_0,0)$, which causes the singularity of the source nonlinearity. Our approach  is then based on the nonexistence for the nonhomogeneous Hardy problem
$$
\cL_\mu^s u\geq g \quad
    {\rm in}\ \ \cD'(\Omega\setminus\{0\}), \qquad
u\geq0 \quad   {\rm in}\ \  \R^N\setminus  \Omega,
$$
when the nonnegative function $g$ is locally H\"older continuous in $\bar\Omega\setminus\{0\}$,
but it doesn't belong to $ L^1(\Omega,\, \Gamma_{s,\mu} dx)$.
More precisely, the main idea can be summarized as follows: We will argue by contradiction,  if (\ref{eq 1.1-in})   admits a positive  solution $u_0$ for $\mu\in[\mu_0,0)$,   which   blows up at least like $\Gamma_{s,\mu}$  at the origin (here the other possibilities are $-\frac{2s+\theta}{p-1}$ or $\Phi_{s,\mu}$).
 Let
\begin{equation}\label{q exp}
  q^\#_{\mu,\theta}:=\frac{N+\theta}{-\tau_+}-1\ \ {\rm if }\ \, \mu\not=0,\qquad
  q^\#_{\mu,\theta}=+\infty\ \ {\rm if }\ \,  \mu=0.
 \end{equation}
 We remark here that for $N+\theta>0$,
 $$q^\#_{\mu,\theta}=p^\#_{\mu,\theta}\ \ {\rm if}\ \, \mu=\mu_0,\quad\  q^\#_{\mu,\theta}>p^\#_{\mu,\theta}\ \ {\rm if}\ \, \mu_0<\mu<0 \quad\  q^\#_{\mu,\theta}<0\ \ {\rm if} \ \, \mu>0.$$

Then for $p\geq q^\#_{\mu,\theta}$,  a contradiction follows thanks to the fact that $Qu_0^p\not\in L^1(\Omega,\, \Gamma_{s,\mu} dx)$.
 For $\mu\in(\mu_0,0)$ and  $p\in(p^\#_{\mu,\theta}, q^\#_{\mu,\theta})$, we  improve the blowing up rate up to $\tau_-$ iteratively, then a contradiction comes from the fact that $Qu_0^p\not\in L^1(\Omega,\, \Gamma_{s,\mu} dx)$. For the critical case $p=p^\#_{\mu,\theta}>1$, we use the observations:
$\displaystyle \liminf_{|x|\to0}Q(x)u_0^{p-1}(x)|x|^{2s}>0$ we then use that, $ \frac{1}{2}Qu_0^{p-1}$ decreases  the coefficient $\mu$ and  then the problem (\ref{eq 1.1-in}) with $\mu$ is turned to the same one replaced by $\mu'(<\mu)$, then  $p=p^\#_{\mu,\theta}$ is supercritical, which implies the nonexistence.  \smallskip

However, for the critical case $p=p^\#_{\mu,\theta}\in(0,1)$ and $\mu\in(\mu_0,0)$, it is a totally different phenomena:
\begin{theorem}\label{teo cri ex 1}
Assume that  $\mu\in(\mu_0, 0)$,  (\ref{pt r1})  holds for some  $\theta\in\R$
and
\begin{equation}\label{pt r1-l}
  \limsup_{|x|\to0^+}Q(x)|x|^{-\theta}<+\infty.
\end{equation}
%Let $p^\#_{\theta,\mu}$ be defined in (\ref{eq 1.1-cr1}).\smallskip

  If $\theta\in(\tau_+-2s, -2s)$,
  then $p^\#_{\mu,\theta}\in(0,1)$ and $p=p^\#_{\mu,\theta}$,
   then problem (\ref{eq 1.1-in})    has infinitely many  positive classical solutions.

\end{theorem}

For $\mu\geq 0$, we have the following nonexistence result in a punctured domain:

 \begin{theorem}\label{teo 1-1}
Assume that  $\mu\geq  0$  and $Q$ verify (\ref{pt r1})  with $\theta< -\tau_+-N $, $ q^\#_{\mu,\theta}$ is given by (\ref{q exp}).

 Then for $0<p\leq q^\#_{\mu,\theta}$,    problem (\ref{eq 1.1-in})    has no positive solution.
 \end{theorem}

 \smallskip

Our next interest is to study the nonexistence of (\ref{eq 1.1}) in the exterior domain $\R^N\setminus\bar \Omega$, subject to nonnegative boundary condition, i.e.
\begin{equation}\label{eq 1.1-ext}
\left\{\arraycolsep=1pt
\begin{array}{lll}
 \cL_\mu^s u\geq Q    u^p \quad
   &{\rm in}\ \  \R^N\setminus \bar\Omega \\[2mm]
 \phantom{ \cL_\mu^s   }
u\geq0 \quad &{\rm in}\ \ \,  \bar\Omega.
 \end{array}
 \right.
\end{equation}
Here we use the Serrin's critical exponent $p^*_{\mu,\theta}=1+\frac{2s+\theta}{-\tau_-}$ for $\mu\geq\mu_0$, which has the following properties:
$$p^*_{\mu,\theta}> 1\ \  {\rm if}\ \  \theta> -2s,\quad\  p^*_{\mu,\theta}= 1\ \  {\rm if}\ \  \theta= -2s,\quad\  p^*_{\mu,\theta}\in(0,1)\ \   {\rm if}\ \    \tau_--2s<\theta< -2s.$$

  \begin{theorem}\label{teo 2}
Let    $\mu\geq \mu_0$,   $Q$ verify (\ref{pt r2})  with $\theta>\tau_--2s$ and   $p^*_{\mu,\theta}$ be defined in (\ref{eq 1.1-cr1}).\smallskip

  $(i)$ If $\mu= \mu_0$, $\theta>\tau_--2s$ and  $p\in(0, p^*_{\mu,\theta}]$,  then problem (\ref{eq 1.1-ext})  has no positive weak solution;\smallskip

 $(ii)$ If $\mu>\mu_0$, $\theta>-2s$ and  $p\in(0, p^*_{\mu,\theta}]$,  then problem (\ref{eq 1.1-ext})  has no positive weak solution;\smallskip

$(iii)$ If $\mu>\mu_0$, $\theta= -2s$, either $p<1$ or
 $$p=1\ \ \  {\rm and}\ \ \  \liminf_{|x|\to+\infty}Q(x)|x|^{2s}>\mu-\mu_0,$$
  then problem (\ref{eq 1.1-ext})  has no positive weak solution;\smallskip

 $(iv)$ If $\mu>\mu_0$, $\tau_--2s<\theta< -2s$ and  $p\in(0, p^*_{\mu,\theta})$, then (\ref{eq 1.1-ext})  has no positive weak solution.

\end{theorem}

At the first glance, one possible way for the nonexistence in exterior domain $\R^N\setminus B_1(0)$  is to
use the Kelvin transformation. In fact,  let
$$
  u_0^\sharp(x)=|x|^{2s-N}u_0(\frac{x}{|x|^2})\quad{\rm for}\  x\in \R^N\setminus\{0\}
 $$
then  we have that
$$\cL_{\mu}^s u_0^\sharp(x)=Q(\frac{x}{|x|^2})|x|^{p(N-2s)-2s-N} (u_0^\sharp(x))^p   \quad{\rm in}\ \ B_1\setminus\{0\}.$$
Since the potential $Q(\frac{x}{|x|^2})|x|^{p(N-2s)-2s-N}\geq c|x|^{p(N-2s)-2s-N-\theta}$ however this does not seem to be conclusive as the exponent $p(N-2s)-2s-N-\theta$ becomes complicated compared with $-2s$. For this, our idea for the exterior domain is to obtain an initial decay $|x|^{\tau_-}$ at infinity, a contradiction will be obtained with the nonexistence of nonhomogenous problem in the interior domain, if not, due to the interaction between the fractional Hardy operator and the nonlinearity,   we  can improve  the decay  at infinity until  $Qu_0^p$ with the new decay  which does not correspond to the requirements of the nonhomogeneous problem in the exterior domain.

For $p=p^*_{\theta,\mu}\in(0,1)$, we have the following existence result:
 \begin{theorem}\label{teo cri ex 2}
Assume that     $\mu>\mu_0$ and there is  $\theta\in\R$ such that  (\ref{pt r2}) holds and
\begin{equation}\label{pt r1-l}
  \limsup_{|x|\to0^+}Q(x)|x|^{-\theta}<+\infty.
\end{equation}
%Let  $p^*_{\theta,\mu}$ be defined in (\ref{eq 1.1-cr1}).\smallskip

  If $\tau_--2s<\theta< -2s$ and  $p= p^*_{\mu,\theta}\in(0,1)$, then (\ref{eq 1.1-ext})  has  infinitely many  positive classical solutions.

\end{theorem}

Our methods allow us  to deal with the case: $p\in(0,1)$ and $\theta<-2s$.
 To the best of our knowledge, our critical exponents are sharp and our method are new and very robust, especially for $\theta<-2s$.  Our Liouville type theorems can be applied to classify the positive solution of the fractional Hardy problem in $\R^N\setminus\{0\}$, i.e.
  \begin{equation}\label{eq 1.1 whole}
\cL_{\mu}^s u=    |x|^{\theta} u^p\quad     {\rm in}\ \,  \R^N\setminus\{0\}.
\end{equation}

\begin{corollary} \label{cr 1.1}

$(i)$ When $\mu=\mu_0$,
 problem (\ref{eq 1.1 whole}) has no positive weak solution   for any $p\in(0,+\infty)$ and any $\theta\in\R$.
 \smallskip

$(ii)$ {\it Non-existence:}  When $\mu>\mu_0$, problem (\ref{eq 1.1 whole}) has no positive weak solution if  \\[1mm] Either:   {\it Case 1: $\mu\in(\mu_0,+\infty)$ }
$$  p\in(0,p^*_{\mu,\theta} ] \ \ {for}\ \ \theta> -2s,\qquad p\in(0,p^*_{\mu,\theta} ) \quad{for}\ \ \tau_--2s<\theta< -2s;$$
or {\it Case 2:   $\mu\in(\mu_0,0)$ }
$$p\in [p^\#_{\mu,\theta} ,+\infty)\ \, {\rm for \ } \,  \theta> -2s,\quad\, p\in (p^\#_{\mu,\theta} ,+\infty)\ \, {\rm for \ } \, \theta\in(\tau_+-2s,-2s],$$
$$  p>0\ \, {\rm for \ }\,  \theta\leq \tau_+-2s;$$

 $(iii)$ {\it Existence:} When $\mu>\mu_0$
 problem (\ref{eq 1.1 whole}) has infinitely many   positive classical solutions if one of the following holds:  \smallskip

  \ \ \ \ \ \quad {\rm (a)}\  $\theta\geq \tau_--2s$ and $p\in (p^{*}_{\mu,\theta} ,\, p^\#_{\mu,\theta} )\cap(0,+\infty)$;  \smallskip

   \ \ \ \ \ \quad {\rm (b)}\ $\mu\in(\mu_0,0)$, $\theta\in (\tau_+-2s,\, \tau_+-2s)$ and $p=p^\#_{\mu,\theta}$; \smallskip

   \ \ \ \ \ \quad {\rm (c)}\ $\mu> \mu_0$, $\theta\in (\tau_--2s,\, \tau_--2s)$ and $p=p^*_{\mu,\theta}$.
\end{corollary}

\begin{remark}
$(i)$ In a punctured domain, \cite{F} investigated the nonexistence of (\ref{eq 1.1-in}) for $\mu\in[\mu_0,0)$ and $Q=1$,  by using  the extension of the fractional Laplacian to a local problem in the half space by adding a dimension.  In exterior domains, \cite{W} obtained nonexistence of solutions of (\ref{eq 1.1-ext}) in the  case: $\theta>-2s$ and  $p\in(0, p^*_{\mu,\theta})$.

Our results also hold for $s=1$. \smallskip

\noindent $(ii)$ It is an interesting phenomena for the existence of super solution in the case: $\tau_+-2s<\theta<-2s$, $p^\#_{\mu,\theta}\in(0,1)$ and in the one: $\theta\in (\tau_--2s,\, \tau_--2s)$ and $p=p^*_{\mu,\theta}$.
\smallskip

\noindent
%$(iii)$ When $\theta>\tau_+-2s$, another interesting model is
%$$\cL^s_\mu u=\big(1+\big|\ln |x| \big|\big)^{\tau}   |x|^{\theta} u^{p^\#_{\theta,\mu}} $$ in a punctured domain or in a exterior domain, what is the critical exponent of $\tau$ for the nonexistence.

\end{remark}

The rest of this paper is organized as follows. In Section 2, we classify the nonexistence of positive super solutions for  nonhomogeneous fractional Hardy problems, provide necessary estimates and prove Corollary \ref{cr 1.1}.    Section 3 and Section 4 are devoted to prove the nonexistence of  positive super solutions in bounded punctured domain and  in an exterior domain respectively, i.e. the proofs of Theorem
\ref{teo 1} and Theorem \ref{teo 2}. Finally, we show the existence in the critical case in Theorem \ref{teo cri ex 1} and Theorem \ref{teo cri ex 2}.

\setcounter{equation}{0}
\section{ Preliminary  }

\subsection{Basic results  }
We start this section by recalling the classification of nonhomogeneous Hardy problems
\begin{equation}\label{eq 2.1fk}
 \cL^\s_\mu u  = f   \ \ {\rm in}\quad \Omega\setminus\{0\},\qquad
u=0 \ \ {\rm in}\ \ \R^N \setminus \Omega.
\end{equation}
 Let us introduce a potential function
 $$
 \Lambda _\mu(x)= \left\{\arraycolsep=1pt
\begin{array}{lll}
  \displaystyle 1\quad
   &{\rm if}\ \,  \tau_+ >2\s-1,\\[1mm]
   |x|^{1-2\s+\tau_+ } \quad
   &{\rm if}\ \,  \tau_+<2\s-1,\\[1mm]
 \phantom{   }
 \displaystyle  1+(-\ln|x|)_+ \quad &{\rm   if } \ \,  \tau_+=2\s-1
 \end{array}
 \right.
$$
and the following theorem follows \cite[Theorem 1.3]{CW} directly.

\begin{theorem}
\label{theorem-C}
Let $\mu\ge\mu_0$ and $f\in C^\beta_{loc}(\bar\Omega\setminus \{0\})$ for some $\beta\in(0,1)$.
\begin{enumerate}
\item[(i)] (Existence) If $f  \in L^1(\Omega,\Gamma_{s,\mu}(x) dx)$, then for every $k\in\R$ there exists a solution $u_k \in L^1(\Omega,\Lambda_\mu  dx)$ of   problem (\ref{eq 2.1fk})
satisfying the distributional identity
\begin{equation}
  \label{eq:distributional-k}
\int_{\Omega}u_k   (-\Delta)^\s_{\Gamma_{s,\mu}} \xi \,dx = \int_{\Omega}f \xi\, \Gamma_{s,\mu} dx +c_{\s,\mu} k\xi(0) \qquad  {\rm for\ all}\ \,\xi \in  C^2_0(\Omega).
\end{equation}
\item[(ii)]  (Existence and Uniqueness) If $f \in L^\infty(\Omega, |x|^{\rho}dx)$ for some $\rho < 2\s - \tau_+(\s,\mu)$, then for every $k\in\R$ there exists a unique solution $u_k \in L^1(\Omega,\Lambda_\mu  dx)$ of  problem (\ref{eq 2.1fk}) with the asymptotics
 \begin{equation}\label{beh 1}
 \lim_{|x| \to 0^+}\:\frac{u_k(x)}{\Phi_{s,\mu}(x)} = k.
 \end{equation}
Moreover, $u_k$ satisfies the distributional identity (\ref{eq:distributional-k}).
\item[(iii)] (Nonexistence) If $f$ is nonnegative and
\begin{equation}\label{f2}
 \lim_{\epsilon\to0^+}\int_{\Omega\setminus B_\epsilon(0)} f\, \Gamma_{s,\mu} dx =+\infty,
\end{equation}
then the problem
\begin{equation}\label{eq 1.1f}
 \arraycolsep=1pt\left\{
\begin{array}{lll}
 \displaystyle  \mathcal{L}_\mu^\s  u= f\quad
   {\rm in}\ \, {\Omega}\setminus \{0\},\\[1.5mm]
 \phantom{  L_\mu \, }
 \displaystyle  u\ge 0\quad  {\rm   in}\ \,  \R^N\setminus \Omega
 \end{array}\right.
\end{equation}
has no nonnegative   solution $u \in L^\infty_{loc}(\R^N \setminus \{0\}) \cap L^1_s(\R^N)$.
\end{enumerate}
\end{theorem}

\begin{proposition}\cite[Proposition 2.4]{CV0} \label{Kato 2.1}
If $f\in L^1(\cK,\rho_{_\cK}^s dx)$, there exists a  unique weak
solution $u$ of the problem
\begin{equation}\label{homo}
 (-\Delta)^s  u=f\quad   \rm{in}\quad \cK, \qquad
 u=0\quad   \rm{a.e.\  in}\quad \R^N\setminus\cK.
\end{equation}
For any $\xi\in\mathbb{X}_{s,+}(\cK)$,  we have
 \begin{equation}\label{sign}
\int_\cK |u|(-\Delta)^s \xi dx\le \int_\cK  \xi\
{\rm{sgn}}(u)f \ dx
\end{equation}
and
 \begin{equation}\label{sign+}
\int_\cK u_+(-\Delta)^s \xi dx\le \int_\cK  \xi \
{\rm{sgn}}_+(u)f \ dx.
\end{equation}
where ${\rm{sgn}}$ is the signum  and $t_+=\max\{0, t\}$.
\end{proposition}
Motivated by the above proposition,  we  obtain the following comparison principle in the weak sense.
\begin{lemma}  \label{Com 2.1}
Let $\mu\geq\mu_0$, $\cK$ be a bounded $C^2$ domain such that  $0\not\in \overline{\cK}$ and $u$  be a weak solution of
\begin{equation}\label{homo}
\cL^s_\mu  u\leq 0\quad   \rm{in}\ \ \cK, \qquad
 u\leq 0\quad   \rm{a.e.\ in}\ \ \R^N\setminus\cK
\end{equation}
in the sense that $v\in
L^1_s(\R^N)$ and
$$%\begin{equation}\label{weak sense}
\int_{\cK } v \cL_\mu^s\xi \, dx\leq \int_{\cK } g \xi
dx,\quad \forall\xi\in \mathbb{X}_{s,+}(\cK).
$$%\end{equation}

 Then $u\leq0$ a.e. in $\R^N$.
\end{lemma}
\noindent{\bf Proof. } Let $w$ be a function such that $w=0$ in $\cK$ and $w=u$ in $\cK^c$.
Then for $x\in \cK$, we have that
$$f_0(x):=(-\Delta)^s w(x)=-C_{N,s} \int_{\cK^c} \frac{u(y)}{|x-y|^{N-2s}}dy\leq 0, $$
which is continuous locally in $\cK$.
Let
$$v=u-w$$
then $v$ is a weak solution of
\begin{equation}\label{homo}
(-\Delta)^s  v(x)=-\frac{\mu}{|x|^{2s}}v(x)+f_0(x)  \quad   \rm{in}\ \ \cK, \qquad
 u= 0\quad   \rm{a.e.\ in}\ \ \R^N\setminus\cK.
\end{equation}
From (\ref{sign+}), we obtain that
 \begin{eqnarray*}
\int_\cK v_+(-\Delta)^s \xi dx\le \int_\cK  \xi \
{\rm{sgn}}_+(v) \Big(-\frac{\mu}{|x|^{2s}}v(x)+f_0(x)\Big) \ dx.
\end{eqnarray*}
which implies that
 \begin{eqnarray} \label{kato 1-0}
\int_\cK v_+\cL^s_\mu \xi dx\le \int_\cK  \xi \
{\rm{sgn}}_+(v)  f_0(x)  \ dx.
\end{eqnarray}
For $\mu\geq\mu_0$,  let $\xi_1$ be the unique solution of
 \begin{eqnarray*}
 \cL^s_\mu \xi_1=1\quad{\rm in}\ \, \cK,\qquad \xi_1=0 \quad{\rm in}\ \, \R^N\setminus\cK,
\end{eqnarray*}
which is positive, regular in $\cK$ since $0\not\in \overline{\cK}$.
Then $\xi_1\in \mathbb{X}_{s,+}(\cK)$ and take $\xi=\xi_1$ in (\ref{kato 1-0}) to obtain that
 \begin{eqnarray*}
\int_\cK v_+  dx\le \int_\cK  \xi_1 \
{\rm{sgn}}_+(v)  f_0(x)  \ dx\leq 0.
\end{eqnarray*}
Thus, $u\leq v_+= 0$ a.e. in $\cK$.\hfill$\Box$\smallskip

From above comparison principle in the weak sense, we now consider  the nonexistence of the super weak solution to  \begin{equation}\label{eq 2.2-0}
\cL_{\mu}^s v\geq  g\ \      {\rm in}\ \,   O\setminus\{0\}, \qquad
   v\geq 0\ \    {\rm a.e.\ in}\ \,   \R^N\setminus O,
\end{equation}
where $O$ is a bounded $C^2$ domain containing the origin.
{\it Here $v$ is a weak solution of (\ref{eq 2.2-0})  if $v\in
L^1_s(\R^N)$,  given any  bounded Lipschitz open set $\cK_\epsilon$ defined in (\ref{kkk}) such that
$$%\begin{equation}\label{weak sense}
\int_{\cK_\epsilon} v \cL_\mu^s\xi \, dx\geq \int_{\cK_\epsilon} g \xi
dx,\quad \forall\xi\in \mathbb{X}_{s,+}(\cK_\epsilon).
$$%\end{equation}
The weak sub-solution has the similar definition. }

 In this weak sense, we have the following comparison principle.

\begin{lemma}\label{cr hp}
Assume that $\mu\geq\mu_0$, $O$ is a bounded $C^2$ domain containing the origin  and
\begin{equation}\label{eq0 2.1}
  \mathcal{L}_\mu^s u  \leq 0 \ \
   {\rm in}\quad  O\setminus \{0\},\qquad   u \leq  0\ \   {\rm   in}\ \ \R^N\setminus  O\quad {\rm in\ the \ weak\ sense}
\end{equation}
and
$$\limsup_{|x|\to0^+}u(x)\Phi_{s,\mu}^{-1}(x)\leq 0.$$
 Then
$$u\leq 0\quad{\rm a.e.\ in}\quad O.$$
\end{lemma}
{\bf Proof.}   {\it Case: $\mu>\mu_0$. } $$\limsup_{|x|\to0^+}u(x)\Phi_{s,\mu}^{-1}(x)\leq0, $$
then for any $\epsilon>0$, there exists $r_\epsilon>0$ converging to zero as $\epsilon\to0$ such that
 $$u\le \epsilon \Phi_{s,\mu}\quad{\rm a.e.\ in}\quad \overline{B_{r_\epsilon}}\setminus\{0\}.$$

   Let $v=u-\epsilon \Phi_{s,\mu}$ satisfy that
$$\mathcal{L}_\mu^s v\le 0\quad {\rm  in}\ \, O\setminus \overline{B_{r_\epsilon}} \quad {\rm in\ the\ weak\ sense }  $$
and
$$u=-\epsilon \Phi_{s,\mu}<0 \quad{\rm a.e.\ in}\ \ \R^N\setminus (O\setminus \overline{B_{r_\epsilon}}).$$
Then by   Lemma \ref{Com 2.1},  for $\mu\geq\mu_0$  we have that
 $u\le \epsilon \Phi_\mu $ a.e. in $O\setminus\{0\}. $
Since $\epsilon$ is arbitrary, we have that $u\le 0$ a.e.  in $O$.\smallskip

 {\it Case: $\mu=\mu_0$. }  It is shown \cite[Theorem 4.2]{CW} that
$$\left\{\arraycolsep=1pt
\begin{array}{lll}
\mathcal{L}_\mu^\s  u= 0\quad
   {\rm in}\ \,  O\setminus \{0\},\\[2mm]
\quad \, u= 0\quad  {\rm   in}\ \,  \R^N\setminus O,\\[2mm]
 \displaystyle\lim_{|x|\to0^+}\frac{u(x)}{\Phi_{s,\mu_0}(x)}=1
  \end{array}
 \right.
$$
has a  unique solution $ \Phi_{s,\mu_0}^O$, which is positive in $O\setminus\{0\}$. Repeating the above argument of the case $\mu>\mu_0$ replacing
$\Phi_{s,\mu}$ by $\Phi_{s,\mu_0}^O$.
  \hfill$\Box$

  \begin{remark}
  The above corollary holds in the classical sense and in this case we just use the comparison
  principle by applying the fractional Hardy inequality \cite{BM}.

  \end{remark}

 \begin{theorem}\label{teo 2.1} Assume that $g $ is a nonnegative function in $L^1_{loc}(\bar O\setminus\{0\})$.

Then  problem (\ref{eq 2.2-0})   has no positive weak solutions   if
 $$\lim_{r\to0^+}\int_{B_{r_0}(0)\setminus B_r(0)}g(x)\Gamma_{s,\mu} dx=+\infty,$$
 where $r_0>0$ such that $B_{r_0}\subset O$.

In particular, the above assumption can be replaced by
$$\liminf_{|x|\to0^+} g(x)|x|^{ 2s-\tau_-} >0.$$

\end{theorem}
\noindent{\bf Proof. } By contradiction, we assume that $u$ is a positive solution of (\ref{eq 2.1-hom}). Let $g_0$ be a function in  $C^\beta_{loc}(O\setminus \{0\})$ for some $\beta\in(0,1)$  such that $0\leq g_0\leq g$ a.e. in $O$ and
$$\lim_{r\to0^+}\int_{B_{r_0}\setminus B_r}g_0(x)\Gamma_{s,\mu} dx=+\infty.$$
Let $$g_n=g_0\eta_0(n|x|),$$
where $\eta_0:[0,+\infty)\to[0,1]$ be a nondecreasing smooth function such that
\begin{equation}\label{eta}
\eta_0=1\quad{in}\ \ [2,+\infty], \qquad \eta_0=0\quad{\rm in} \ \ [0,1].\end{equation}
Since $g_n$ is H\"older continuous, it follows by \cite[Theorem 4.1]{CW} that
 \begin{equation}\label{eq 2.1-n}
\left\{\arraycolsep=1pt
\begin{array}{lll}
\cL_{\mu}^s v=  g_n\qquad      {\rm in}\ \,  O\setminus\{0\}, \\[2mm]
\quad \,  v=0\qquad   \,  \ {\rm in}\ \,   \R^N\setminus O,\\[2mm]
\displaystyle \lim_{|x|\to0^+}v(x)\Phi_{s,\mu}^{-1}(x)=0
  \end{array}
 \right.
\end{equation}
  has unique solution $v_n$ and by Lemma \ref{cr hp}, it implies that
  $$0<v_{n}\leq    u\quad {\rm in}\ \, O\setminus\{0\}. $$
Standard local regularity results show that  the limit $\{v_n\}_n$, $v_\infty$, is a positive classical solution of
  \begin{equation}\label{eq 2.1-l}
\left\{\arraycolsep=1pt
\begin{array}{lll}
\cL_{\mu}^s v=  g_0\quad \     {\rm in}\ \  O\setminus\{0\},\\[2mm]
\quad\  v=0\quad   \  {\rm in}\ \  \R^N\setminus O,
  \end{array}
 \right.
\end{equation}
then we obtain a contradiction from  Theorem \ref{theorem-C} part $(iii)$.
  \hfill$\Box$\medskip

  \begin{corollary}\label{re 2.1} Let $\mu\ge\mu_0$,  $O$ is a bounded domain containing the origin and nonnegative function  $f\in C^\beta_{loc}(\R^N \setminus \bar O)$ for some $\beta\in(0,1)$.
The homogeneous problem
\begin{equation}\label{eq 1.1 EH}
\left\{\arraycolsep=1pt
\begin{array}{lll}
 \cL^\s_\mu u  \geq  f   \quad
   &{\rm in}\ \  \R^N\setminus\bar O, \\[2mm]
 \phantom{ \cL_\mu^s   }
u\geq 0 \quad &{\rm a.e.\ in}\ \ \,  \bar O
 \end{array}
 \right.
\end{equation}
 has no positive weak solutions if
 $$\lim_{r\to+\infty}\int_{B_r \setminus B_{R_0} }f(x)|x|^{ 2s-\tau_+-N}dx=+\infty.$$
In particular, the above assumption can be replaced by
$$\liminf_{|x|\to+\infty} f(x)|x|^{2s -\tau_+} >0.$$

\end{corollary}
   \noindent{\bf Proof. }  Let $f_0$ be a function in  $ C^\beta_{loc}(\R^N \setminus \bar O)$ for some $\beta\in(0,1)$    such that $0\leq f_0\leq f$ a.e. in $\R^N\setminus B_{R_0}$ and
 $$\lim_{r\to+\infty}\int_{B_r \setminus B_{R_0} }f_0(x)|x|^{ 2s-\tau_+-N}dx=+\infty.$$

By contradiction, we assume that (\ref{eq 1.1 EH}) has a positive weak solution $u_{f}$.

 Let $u_{f_0, n}$ be a positive classical solution of
 $$ \cL^\s_\mu u  =  f_0   \ \  {\rm in}\ \,  B_n \setminus \overline{B_{R_0}}, \qquad u= 0 \ \  {\rm \ in}\ \ \, \overline{B_{R_0}}\cap (\R^N\setminus B_n),$$
 where $n>R_0$.
 We  deduce that
the map $n\to u_{f_0, n}$ is increasing and bounded by $u_{f}$, and from the stability, we have that
\begin{equation}\label{eq 1.1 EH=}
\left\{\arraycolsep=1pt
\begin{array}{lll}
 \cL^\s_\mu u  =  f_0   \quad
   &{\rm in}\ \  \R^N\setminus\bar O, \\[2mm]
 \phantom{ \cL_\mu^s   }
u= 0 \quad &{\rm   in}\ \ \,  \bar O
 \end{array}
 \right.
\end{equation}
has a classical solution $u_{f_0}$. The regularity of $f_0$ and Lemma \ref{cr hp} imply that
$$0\leq u_{f_0}\leq u_f\quad {\rm a.e.\ in}\ \ \R^N.$$

Now denote by
 \begin{eqnarray}\label{kt-1-1}
O^\sharp=\left\{x\in\R^N: \ \frac{x}{|x|^2}\in O\right\}\quad{\rm and}\quad  u^\sharp(x)=|x|^{2s-N}u_{f_0}(\frac{x}{|x|^2})\quad{\rm for}\  x\in O^\sharp.
\end{eqnarray}
  Clearly, $(\R^{N}\setminus\{0\})^\sharp=\R^{N}\setminus\{0\}$  and by the Kelvin transformation \cite[Lemma A.1.1]{CLM}  we obtain that
 \begin{eqnarray}\label{kt-1}
 (-\Delta)^s u^\sharp(x)    =  |x|^{-2s-N} \big((-\Delta)^s u_{f_0}\big)\left(\frac{x}{|x|^2}\right)\quad{\rm for}\  x\in O^\sharp.
 \end{eqnarray}
Therefore,  $u^\sharp$ is a super solution of
\begin{equation}\label{eq 1.1 KT}
\cL_{\mu}^s u^\sharp(x) \geq   |x|^{-2s-N} \tilde {f_0} (x)  \quad    {\rm in}\ \  O^\sharp,
\end{equation}
where $\tilde  f_0 (x)=f_0(\frac{x}{|x|^2})$.

 Note that  there exists $r_1>0$ such that
$$B_{r_1}\subset O^\sharp\cup\{0\}.$$
Note that for $r>\frac1{r_1}>0$
 \begin{eqnarray*}
\int_{B_{r_1}\setminus B_\frac1r }f(\frac{x}{|x|^2})|x|^{ -2s-N+\tau_+}dx&=& \int_{B_r\setminus B_{\frac1{r_1}}}f(y)|y|^{ 2s-\tau_+-N}dy
%\\&=&\int_{B_r(0)\setminus B_{\frac1{r_1}}(0)}f(x)|x|^{ N+\tau_-}dx
\\&\to&+\infty\quad{\rm as}\ \,  r\to+\infty.
 \end{eqnarray*}
The contradiction follows by   Theorem \ref{theorem-C} part $(iii)$.
\hfill$\Box$

 \begin{corollary}\label{re 2.2} Let  $O$ be a bounded domain containing the origin and   $f\in L^1_{loc}(\R^N \setminus \bar O)$  be a nonnegative nontrivial function.

If $\mu<\mu_0$, then the homogeneous problem
\begin{equation}\label{eq 1.1 EH-1}
 \cL^\s_\mu u  \geq  f   \ \ {\rm in}\ \ \R^N\setminus\bar O, \qquad
  u\geq 0 \ \  {\rm in}\ \ \,  \bar O
\end{equation}
 has no positive weak solution.
 \end{corollary}
 {\bf Proof. } By contradiction,  we assume that  (\ref{eq 1.1 EH-1}) has  a positive weak solution $u_f$.
  Let $f_0$ be a function in  $ C^\beta_{loc}(\R^N \setminus \bar O)$ for some $\beta\in(0,1)$    such that $0\leq f_0\leq f$ a.e. in $\R^N\setminus B_{R_0}$
 From the proof of Corollary \ref{re 2.1}, we have that
 \begin{equation}\label{eq 1.1 EH-2}
 \cL^\s_\mu u  =  f_0   \ \ {\rm in}\ \ \R^N\setminus\bar O, \qquad
  u= 0 \ \  {\rm in}\ \ \,  \bar O
\end{equation}
  has a classical solution $u_{f_0}$.

 Recall that
  $u^\sharp$, defined in (\ref{kt-1-1}),  is a super solution of
$$
\cL_{\mu}^s u^\sharp(x) =  |x|^{-2s-N}   f_0 (\frac{x}{|x|^2})  \quad    {\rm in}\ \  (\R^N\setminus\bar O)^\sharp,
$$
where $(\R^N\setminus\bar O)^\sharp$ is a bounded domain containing $B_r\setminus\{0\}$ for some $r>0$.
A contradiction comes from  \cite[Theorem 1.4]{CW}, that says that
 the above problem has no positive classical solution for  $\mu<\mu_0.$ The proof is now complete.\hfill$\Box$

 \subsection{Asymptotics}

In order to improve the blowing up rate at the origin or decay at infinity, we need the following estimates.  From \cite[Lemma 1.1]{CW}
letting
\begin{equation}
  \label{eq:fractional-power}
 c_\s(\tau) = 2^{2\s} \frac{\Gamma(\frac{N+\tau}{2})\Gamma(\frac{2\s-\tau}{2})}{\Gamma(-\frac{\tau}{2})\Gamma(\frac{N-2\s+\tau}{2})}>0,
\end{equation}
the map: $\tau\mapsto c_\s(\tau)$ is concave and
$\tau_+\geq \tau_-$ are the two zero points of
$$c_\s(\tau)-\mu=0 \quad{\rm for}\ \ \mu>\mu_0.$$
There holds that
\begin{equation}
  \label{eq:fractional-power-1}
 (-\Delta)^\s |\cdot|^\tau = c_\s(\tau)|\cdot|^{\tau-2\s} \qquad  {\rm in}\ \  \mathcal{S}'(\R^N).
\end{equation}
 For $\tau\in (\tau_-, \tau_+)$, we set $b_s(\tau):=c_\s(\tau)-\mu>0$.

   \begin{lemma}\label{teo 2.2} Assume that $\mu>\mu_0$,
   $O_r=B_r \setminus \{0\}$ or $O_r=\R^N\setminus\bar B_\frac1r$, the nonnegative function $g\in C^\beta_{loc}(O_r)$ for some $\beta\in(0,1),$ and there exists $\tau\in(\tau_- ,\tau_+ )$, $c>0$ and $r_0>0$ such that
 $$ g(x)\geq c|x|^{\tau-2s}\quad{\rm in}\ \, O_r.$$

Let  $u_g$ be a positive weak solution of
 \begin{equation}\label{eq 2.1-hom}
\cL_\mu^s u\geq g \quad {\rm in}\ \ O_r,\qquad  u\geq0 \quad  {\rm in}\ \  \R^N\setminus   O_r,
 \end{equation}
then there exists $c_1>0$  such that
$$u_g(x)\geq c_1|x|^\tau\quad{\rm in}\ \  O_{\frac r2}.$$

\end{lemma}
{\bf Proof. }  For $\tau\in(\tau_- ,\tau_+ )$, we have that
$$\cL^s_{\mu} |x|^{\tau}=b_s(\tau) |x|^{\tau-2s}\quad  {\rm in}\ \  \R^N\setminus\{0\},$$
where $b_s(\tau)>0$.

In the case $O=B_1\setminus \{0\}$, we use the function
$$w(x)=|x|^{\tau}-|x|^{\tau_+}\quad  {\rm in}\ \  \R^N\setminus\{0\}$$
as a sub-solution
$$\cL^s_{\mu} u=c_s(\tau)|x|^{\tau-2s}\quad{\rm in}\ B_1\setminus\{0\},\qquad
u\leq 0\ \ {\rm in} \ \ \R^N\setminus B_1,$$
where $c_s(\tau)>0$.
 Then our argument follows by Lemma \ref{cr hp}.

 When $O=\R^N\setminus B_1$,  we use the function
$$w(x)=|x|^{\tau}-|x|^{\tau_-}\quad  {\rm in}\ \  \R^N\setminus\{0\}$$
as a sub-solution
$$\cL^s_{\mu} u=c_s(\tau)|x|^{\tau-2s}\quad{\rm in}\ \R^N\setminus \bar B_1,\qquad
u\leq 0\ \ {\rm in} \ \  \bar B_1$$
 and the rest of the proof is standard.  \hfill$\Box$

 \subsection{Liouville Theorem in $\R^N\setminus\{0\}$.  }

\noindent {\bf Proof of Corollary \ref{cr 1.1}. }    For $\mu\geq \mu_0$, we can see that
$$ p^*_{\mu,\theta}\leq \frac{N+2s+2s\theta}{N-2s} \leq p^\#_{\mu,\theta},$$
where the equality hold only when $\mu=\mu_0$.

We first deal with the case: $\mu>\mu_0$.\smallskip

{\it Nonexistence: } Thanks to the domain $\R^N\setminus\{0\}$ and $Q(x)=|x|^{\theta}$,    the nonexistence of positive solutions to
 $$
\cL_{\mu}^s u=    |x|^{\theta} u^p\quad     {\rm in}\ \,  \R^N\setminus\{0\}
$$
  follows from Theorem \ref{teo 2} for
  $$p\in(0, p^*_{\mu,\theta}]\quad {\rm if} \ \ \theta>-2s$$
  or for
   $$ p\in(0, p^*_{\mu,\theta})\quad {\rm if} \ \ \theta\in(\tau_+-2s,-2s).$$
By Theorem \ref{teo 1},  (\ref{eq 1.1 whole}) has no positive weak solutions if
$$p\in [p^\#_{\mu,\theta} ,+\infty)\quad{\rm for \ } \mu\in(\mu_0,0),\ \theta> -2s$$
or
$$p\in (p^\#_{\mu,\theta} ,+\infty)\cap (0,+\infty)\quad{\rm for \ } \mu\in(\mu_0,0),\ \theta\leq -2s.$$

 \smallskip

{\it Existence:  Part $(a)$. }  Note that for $\tau\in (\tau_-,\tau_+)$,
$$b_s(\tau):= c_\s(\tau)-\mu>0.$$
Note that when  $p\in(p^*_{\mu,\theta},\, p^\#_{\mu,\theta})$, we have that
$$-\frac{2s+\theta}{p-1}\in (\tau_-,\tau_+).$$
Let $\tau_p=-\frac{2s+\theta}{p-1}$,   then
$$u_p(x):=c_\s(\tau_p)^{\frac{1}{p-1}} |x|^{\tau_p}$$
 is
a positive classical solution of (\ref{eq 1.1 whole}).   \smallskip

Part $(b)$ and $(c)$ follow from Theorem \ref{teo cri ex 1} and from Theorem \ref{teo cri ex 2} respectively. \smallskip

When $\mu=\mu_0$,
the nonexistence follows by Theorem \ref{teo 1} part $(i)$ and Theorem \ref{teo 2} part $(i)$  directly by the fact that $p^*_{\mu_0,\theta}= p^\#_{\mu_0,\theta}.$
\hfill$\Box$

\setcounter{equation}{0}
\section{ Nonexistence in  a punctured domain }
\subsection{The case $\mu\in[\mu_0,0)$}
We prove the nonexistence of positive solutions of (\ref{eq 1.1-in}) in a punctured domain $\Omega\setminus\{0\}$ by contradiction, i.e.   (\ref{eq 1.1-in}) is assumed to have a positive solution $u_0$ and we will obtain a contradiction from  Theorem \ref{teo 2.1}.
In this section, we   assume that  $$B_4\subset \Omega$$
and
$$Q(x)\geq q_0|x|^{\theta}\quad{\rm in}\ \, B_3\setminus\{0\}.$$

 %  Since $\partial B_{r_0}$ is compact and $u_0$ is regular and positive,
  % then it is well-defined that  $t_0>0$ such that
  % $$t_0=\min_{x\in \partial B_{r_0}(0)} u_0(x).$$

 \begin{proposition}\label{lm 2.1}
Let  $\mu\ge \mu_0$,  $\theta\in\R$ and $u_0$ be a positive weak solution of (\ref{eq 1.1-in}) in $\Omega\setminus\{0\}$,
 then
 $$\liminf_{|x|\to0^+} u_0(x) |x|^{-\tau_+}>0.$$
 \end{proposition}

  In order to obtain the lower bound, we first consider the
one for nonhomogeneous problem
  \begin{equation}\label{eq 2.1-green}
\left\{\arraycolsep=1pt
\begin{array}{lll}
\cL_{\mu}^s v=  f  \qquad      {\rm in}\ \   \Omega\setminus\{0\}, \\[2mm]
\quad \,  v=0\qquad      {\rm in}\ \    \R^N\setminus \Omega.  \end{array}
 \right.
\end{equation}
  {\it Here $v$ is a weak solution of (\ref{eq 2.1-green})  if $v\in
L^1_s(\R^N)$,  given any  bounded Lipschitz open set $\cK$ such that $0\not\in  \overline \cK,\  \overline \cK \subset  O$ and}
$$%\begin{equation}\label{weak sense}
\int_\cK  v \cL_\mu^s\xi \, dx = \int_\cK f \xi
dx,\quad \forall\xi\in \mathbb{X}_{s}(\cK).
$$

 \begin{lemma}\label{lm 2.1-GK}
Let  $\mu\ge \mu_0$, $f\in L^1_{loc}(\Omega\setminus\{0\}) $ be a nonnegative  function,  and $v$ be a weak solution of the problem (\ref{eq 2.1-green}) satisfying:
$$ \liminf_{|x|\to0^+}v(x)\Phi_{s,\mu}^{-1}(x)\geq 0.$$
If there exists $\varepsilon_0>0$ and $r_1>r_2>0$  such that
$$f\geq \varepsilon_0\quad {\rm in}\ \ B_{r_1}\setminus B_{r_2},$$
then there exist $c_2>0$ and  $\sigma>0$ such that
\begin{equation}\label{e 2.2}
v\geq c_2\Gamma_{s,\mu}\quad{\rm for}\ \, x \in B_{\sigma}\setminus\{0\}.
\end{equation}
 \end{lemma}
{\bf Proof. } Thanks to the fact that
$\cL^s_{\mu}\Gamma_{s,\mu}=0$ in $\R^N\setminus\{0\}$, we would like to
construct a nonnegative  function $v_0$ with compact support in $\Omega$ such that
 \begin{equation}\label{e 2.3-0}
v_0=\Gamma_{s,\mu}\quad {\rm in}\ \ B_{r_1} \setminus\{0\}
 \end{equation}
and
 \begin{equation}\label{e 2.3-00}
\cL^s_{\mu}v_0\leq 0 \quad {\rm in}\ \ B_{r_2} \setminus\{0\} \cup (\R^N\setminus  B_{r_1}).
 \end{equation}
Without loss of generality, we set
$$r_1=4,\ \  r_2=1.$$

Now recall that  $\eta_0:\R_+\to [0,1]$ is a nondecreasing smooth function satisfying (\ref{eta}).
Denote
$$v_1(x)=\Gamma_{s,\mu}(x)\eta_0(|x|),\quad\forall\, x\in\R^N\setminus\{0\}.$$
Direct computation shows that for $x\in B_1 \setminus\{0\}$
$$(-\Delta)^s v_1(x)=(-\Delta)^s \Gamma_{s,\mu}(x) +C_{N,s}\int_{\R^N\setminus B_2 }\frac{\Gamma_{s,\mu}(y)(1-\eta_0(y))}{|x-y|^{N+2s}}dy,$$
where
$$%\begin{eqnarray*}
\Big|C_{N,s}\int_{\R^N\setminus B_2 }\frac{\Gamma_{s,\mu}(y)(1-\eta_0(y))}{|x-y|^{N+2s}}dy \Big|
 \leq     C_{N,s}\int_{\R^N\setminus B_2 }\frac{\Gamma_{s,\mu}(y)}{ (|y|+1)^{N+2s}}dy =:C_1.
$$%\end{eqnarray*}
Therefore, we have that
\begin{equation}\label{e 2.3-1}
\cL^s_\mu v_1 \leq C_1\quad {\rm in}\ \  B_1.
\end{equation}
Moroover, for $x\in \R^N\setminus B_3$, we can see that
$$(-\Delta)^s v_1(x)=  C_{N,s}\int_{\R^N\setminus B_2}\frac{-\Gamma_{s,\mu}(y) \eta_0(y) }{|x-y|^{N+2s}}dy  <0.$$

In order to decrease the estimate (\ref{e 2.3-1}), we let $v_2$ be a nonnegative,  radially symmetric, smooth function such that
$$v_2=1\quad{in}\ \ B_{3}(0)\setminus B_2, \qquad v_2=0\quad{\rm in} \ \ B_1\cup  (\R^N\setminus B_4).$$
From the definition of fractional Laplacian,  there exists $C_2>0$
$$\cL_\mu^s v_2=(-\Delta)^s v_2\leq -C_2\quad {\rm in}\ \,  B_1\cup  \big(\Omega \setminus B_4\big). $$

Now let
$$v_0=\frac{v_1}{C_1}+\frac{v_2}{C_2},$$
which verifies (\ref{e 2.3-0}) and   (\ref{e 2.3-00}). By the smoothness of $v_0$, there exists
$C_3>0$ such that
$$\cL_\mu^s v_0\leq C_3\quad {\rm in}\ \ B_4\setminus B_1.$$
By the comparison principle, see Lemma \ref{cr hp},  there exists $c_3>0$ such that
$$v\geq c_3v_0\quad {\rm in}\ \, \Omega\setminus \{0\}$$
and then   (\ref{e 2.2}) follows. \hfill$\Box$\medskip

\noindent{\bf Proof of Proposition \ref{lm 2.1}. }   Let $f_0$ be a nonnegative nontrivial function in $C^1(\Omega)$  such that
$0\leq f_0(x)\leq Q(x)|x|^{\theta} u_0(x)^p$ a.e. in $\Omega$. Then it follows from \cite[Theorem 1.3 (i)]{CW} in the case $k=0$ that
 \begin{equation}\label{eq 2.1}
\left\{\arraycolsep=1pt
\begin{array}{lll}
\cL_{\mu}^s v=  f_0\quad \ \,    {\rm in}\ \, \,  \Omega\setminus\{0\},\\[2mm]
\quad\   v=0\quad  \ \,     {\rm in}\ \, \,  \R^N\setminus \Omega,\\[2mm]
\displaystyle \lim_{|x|\to0^+}v(x)\Phi_{s,\mu}^{-1}(x)=0
  \end{array}
 \right.
\end{equation}
  has unique solution $v_{f_0}\in C(\Omega\setminus\{0\})$, which is positive in $\Omega\setminus\{0\}$.
   By Lemma \ref{cr hp}, we have that $u_0\geq v_{f_0}$ a.e. in $\Omega$.

Now we can affirm that  $f_0$ has a positive low bound in $B_{\frac74 r_0}\setminus B_{\frac54r_0}$,
since $v_{f_0}$ is positive and continuous in $\Omega\setminus\{0\}$ and
$$Q(x)|x|^{\theta} u_0(x)^p\geq Q(x)|x|^{\theta} v_{f_0}(x)^p.$$
Then by (\ref{e 2.2})  we obtain that
   \begin{equation}\label{e 2.3}
   \liminf_{|x|\to0^+}v_{f_0}(x)\Gamma_{s,\mu}^{-1}(x)>0.
   \end{equation}

  Note that $\displaystyle \liminf_{|x|\to0^+}u_0(x)\Phi_{s,\mu}^{-1}(x)\geq 0$, we deduce by Lemma \ref{cr hp} that
    $$u_0\geq v_{f_0}\quad {\rm in}\ \, \Omega\setminus\{0\},$$
which, along with   (\ref{e 2.3}), implies that
 $$\liminf_{|x|\to0^+} u_0(x) |x|^{-\tau_+}>0.$$
 The proof is complete. \hfill$\Box$\bigskip

Let $\tau_0<0$ and   $\{\tau_j\}_j$ be the sequence generated by
\begin{equation}\label{2.1-jj}
 \tau_j=2s+\theta+p\tau_{j-1}\quad{\rm for}\quad  j=1,2,3\cdots.
\end{equation}
\begin{lemma}\label{lm 2.1--j}

$(i)$ Assume that $\theta\in \R$,
$$%\begin{equation}\label{2.00}
 p> 1+\frac{2s+\theta}{-\tau_0} \quad{\rm and}\quad p\geq 1,
$$%\end{equation}
  then $\{\tau_j\}_j$ is a decreasing sequence and for any $\bar \tau<\tau_0$
 there exists $j_0\in\N $  such that
 \begin{equation}\label{2.3-bp}
 \tau_{j_0}\le \bar \tau  \quad {\rm and}\quad \tau_{j_0-1}>\bar \tau.
 \end{equation}

$(ii)$ Moreover, assume that $\theta+2s<0$ and $p\in(0,1)$,   then $\{\tau_j\}_j$ is a decreasing sequence and for any $\bar \tau\in (\frac{2s+\theta}{1-p}, \tau_0)$
there exists $j_0\in\N $  verifying (\ref{2.3-bp}).
 \end{lemma}
 {\bf Proof. }
For $p> \max\Big\{1,1+\frac{2s+\theta}{-\tau_0}\Big\}$, we have that
$$\tau_1-\tau_0=2s+\theta+\tau_0(p-1)<0$$
and
\begin{equation}\label{2.2-in}
\tau_j-\tau_{j-1} = p(\tau_{j-1}-\tau_{j-2})=p^{j-1} (\tau_1-\tau_0),
\end{equation}
which implies that the sequence $\{\tau_j\}_j$ is  decreasing.
Thanks to the fact that $p\geq 1 $, our conclusions are straightforward.  If $p\in(0,1)$,
we can deduces from (\ref{2.2-in}) that
\begin{eqnarray*}
\tau_j  &=&  \frac{1-p^j}{1-p}(\tau_1-\tau_0)+\tau_0\\
&\to&\frac{1}{1-p}(\tau_1-\tau_0)+\tau_0=\frac{2s+\theta}{1-p}\quad {\rm as}\quad j\to+\infty,
\end{eqnarray*}
 then there exists $j_0>0$ satisfying (\ref{2.3-bp}).
\hfill$\Box$\medskip

Now we are in a position to prove Theorem \ref{teo 1}.\medskip

\noindent{\bf Proof of Theorem \ref{teo 1}. }
By contradiction, we assume that $u_0$ is a positive super-weak solution of (\ref{eq 1.1-in}) in $\Omega\setminus\{0\}$.
For $\mu\in[\mu_0,0)$, we have that $\tau_+<0$.

Let
$$q^\#_{\mu,\theta}=  \frac{N+\theta}{-\tau_+}-1.  $$
Then
$$q^\#_{\mu,\theta}=p^\#_{\mu,\theta}\ \ {\rm if}\ \  \mu=\mu_0$$
and
$$q^\#_{\mu,\theta}>p^\#_{\mu,\theta}\ \ {\rm if}\ \  \mu\in(\mu_0,0).$$

Set
$g(x)=Q(x) u_0(x)^{p}$,
then Proposition \ref{lm 2.1} implies that for some $d_0>0$
$$u_0(x)\geq d_0|x|^{\tau_+}\quad{\rm in}\ \ B_{r_0}(0)\setminus\{0\}.$$

{\bf Part 1: $\mu\in[\mu_0,0)$, $\theta\in\R$ and $p\geq q^\#_{\mu,\theta}$, $p>0$. }
Note that
 \begin{eqnarray*}
 \cL^s_\mu u_0(x) \geq  g(x)\geq    d_0^p  |x|^{\theta+\tau_+p}  \quad{\rm in}\ \ B_{r_0}(0)\setminus\{0\},
 \end{eqnarray*}
 where $\theta+\tau_+p+ \tau_+\leq -N$ and
$$\lim_{r\to0^+}\int_{B_{r_0}(0)\setminus B_r(0)} Q(x)|x|^{\tau_+p}\Gamma_{s,\mu} dx=+\infty$$
by the fact that $p\geq q^\#_{\mu,\theta}$.
As a consequence, we see that $u_0$ is a solution of
$$
\cL_\mu^s u_0 \geq  g\quad      {\rm in}\ \,  \Omega\setminus\{0\},\qquad u_0\geq0\quad {\rm in}\ \ \R^N\setminus\{0\}.
$$
 Then there is a contradiction from Theorem \ref{teo 2.1}. \smallskip

 {\bf Part  2: $\mu\in(\mu_0,0)$, $\theta\in\R$ and $p\in \big(p^\#_{\mu,\theta},\, q^\#_{\mu,\theta}\big)$,  $p>0$. }
 For $p\in(0, 1)$ and $p> p^\#_{\mu,\theta}$, we have that
\begin{equation}\label{4.1}
\quad\frac{2s+\theta}{1-p}>\tau_+.
\end{equation}
Let $\tau_0=\tau_+<0$, then
 $$
\cL^s_\mu  u_0(x)  \geq    q_0d_0^p|x|^{p\tau_0+\theta}=q_0d_0^p |x|^{\tau_1-2s}\quad {\rm in}\ \, B_{r_0}\setminus\{0\},
$$
where
$$\tau_1:=p\tau_0+\theta+2s.$$
 By Lemma \ref{teo 2.2}, we have that
$$u_0(x)\geq d_1|x|^{\tau_1}\quad {\rm in}\ \, B_{r_0}\setminus\{0\}.$$
Iteratively, we recall that
$$\tau_j:=p\tau_{j-1}+\theta+2s,\quad j=1,2,\cdots .$$
Note that
$$\tau_1-\tau_0=(p-1)\tau_0+\theta+2s <0$$
for $p\in(p^\#_{\mu,\theta}, q^\#_{\mu,\theta})$.

If $\tau_{1}p   +\theta +2s \leq \tau_-,$
then $\tau_{1}p   +\theta<-N-\tau_+$ and
\begin{eqnarray*}
 \cL^s_\mu u_0(x) \geq  g(x)\geq   q_0 d_1^p  |x|^{\theta+\tau_{1}p}  \quad{\rm in}\ \ B_{r_0}\setminus\{0\}
 \end{eqnarray*}
 and  a contradiction comes from Theorem \ref{teo 2.1}.

If not, we can iterate the above procedure. If
$$\tau_{j+1}:=\tau_jp+\theta+2s\in  (\tau_-,\tau_+),$$
it follows from Lemma \ref{teo 2.2}  that
$$u_0(x)\geq d_{j+1} |x|^{\tau_{j+1}},$$
where
$$\tau_{j+1}=p\tau_j+2s+\theta<\tau_j.$$
If $\tau_{j+1}p+ \tau_++\theta \leq -N,$
then
\begin{eqnarray*}
 \cL^s_\mu u_0(x) \geq    q_0  d_{j+1}^p  |x|^{\theta+\tau_{j+1}p}  \quad{\rm in}\ \ B_{r_0}\setminus\{0\}
 \end{eqnarray*}
 and  a contradiction comes from Theorem \ref{teo 2.1}. Then we are done.

 In fact, this iteration must stop by finite times due to (\ref{4.1}) for $p\in(0,1)$
  and $\tau_j\to-\infty$ as $j\to+\infty$ if $p\geq 1$.
 \smallskip

{\bf Part  3: $\mu\in(\mu_0,0)$, $\theta>-2s$ and $p=p^\#_{\mu,\theta}> 1$. }
Note that
 \begin{eqnarray*}
 \cL^s_\mu u_0(x) \geq  g(x)\geq    q_0 d_0^p  |x|^{\theta+\tau_+p}  \quad{\rm in}\ \ B_{r_0}\setminus\{0\},
 \end{eqnarray*}
 where in this case
 $$\theta+\tau_+p+2s=\tau_+.$$
 For some $\sigma_0>0$
 $$\frac12 Qu_0^{p-1}(x)\geq \sigma_0|x|^{-2s} \quad{\rm in}\ \ B_{r_0}\setminus\{0\}.$$
 Here we can assume that $\mu-\sigma_0\geq \mu_0$, otherwise, we only take a smaller value for $\sigma_0$.
 So we can write problem (\ref{eq 1.1}) as follows:
\begin{equation}\label{eq 1.1 new}
 \cL_{\mu-\sigma_0}^s u_0\geq \frac12 Q(x)   u_0^p \quad {\rm in}\ \  B_{r_0}\setminus\{0\} ,
 \end{equation}
 which is the critical exponent
 $$p^\#_{\mu-\sigma_0,\theta}=  1+\frac{2s+\theta}{-\tau_+(s,\mu-\sigma_0)} <   1+\frac{2s+\theta}{-\tau_+(s,\mu )}=p^\#_{\mu,\theta}, $$
 since $\mu\mapsto \tau_+(s,\mu )$ is decreasing.
Thus it reduces to the case:    {\bf part 2} for (\ref{eq 1.1 new}) for $p=p^\#_{\mu,\theta}$ is supercritical and  a contradiction can then be deduced.\smallskip

{\bf Part  4: $\mu\in[\mu_0,0)$, $\theta=-2s$ and $p=1$. } This reduces to the  Hardy problem
$$(-\Delta)^s u+(\mu+\mu')|x|^{-2s} u =f \quad{\rm in}\ \ B_r\setminus\{0\},$$
where $r>0$,  $f(x)\geq Q(x)-\mu'|x|^{-2s}\geq0$ and $\mu+\mu'<\mu_0$. \cite[Theorem 1.4]{CW} shows that
 the above problem has no positive solution for
 $\mu+\mu'<\mu_0.$
    This completes the proof.  \hfill$\Box$\medskip
 \subsection{The case $\mu\geq 0$}

 \noindent{\bf Proof of Theorem \ref{teo 1-1}. }
By contradiction, we assume that $u_0$ is a positive super-weak solution of (\ref{eq 1.1-in}).
For $\mu>0$, we have that $\tau_+>0$.

For $\mu\geq 0$, $\theta<-\tau_+-N$, observe that
$$   \frac{N+\theta}{-\tau_+}-1>0.   $$

 Since it is known from  Proposition \ref{lm 2.1}  that
$$u_0(x)\geq c_6|x|^{\tau_+}\quad{\rm in}\ \ B_{r_0}(0)\setminus\{0\}$$
for some $c_6>0$,
then we have that
 \begin{eqnarray*}
 \cL^s_\mu u_0(x) \geq     c_6^p  |x|^{\theta+\tau_+p}:=g(x)  \quad{\rm in}\ \ B_{r_0}(0)\setminus\{0\},
 \end{eqnarray*}
 where $\theta+\tau_+p+ \tau_+\leq -N$ and
$$\lim_{r\to0^+}\int_{B_{r_0}(0)\setminus B_r(0)} Q(x)|x|^{\tau_+p}\Gamma_{s,\mu} dx=+\infty$$
by the fact that $p\leq \frac{N+\theta}{-\tau_+}-1$.
As a consequence, we see that $u_0$ is a solution of
$$
\cL_\mu^s u_0 \geq  g\quad      {\rm in}\ \,  \Omega\setminus\{0\},\qquad u_0\geq0\quad {\rm in}\ \ \R^N\setminus\{0\}.
$$
 Then there is a contradiction from Theorem \ref{teo 2.1}. \hfill$\Box$\medskip

 We remark here that the iteration procedure to improve the blow up at the origin
 can't be used in the case when $\mu\geq 0$, due to the fact that $-\tau_+<0$.

 \setcounter{equation}{0}
 \section{Nonexistence in an exterior domain}
  In this subsection, we   deal with semilinear fractional Hardy inequality in $\R^N \setminus \Omega$ and set
   $$\Omega \subset B_1\quad {\rm and}\quad Q(x)\geq q_\infty |x|^{\theta}\quad{\rm in}\ \R^N\setminus B_1.$$

  % Since $\partial B_{R_0}$ is compact and $u_0$ is regular and positive,
 %  then it is well-defined that  $t_0>0$ such that
  % $$t_0=\min_{x\in \partial B_{R_0}(0)} u_0(x).$$

 \begin{proposition}\label{lm 2.1-ex}
Let  $\mu\ge \mu_0$,  $\theta\in\R,$ and $u_0$ be a positive weak solution of (\ref{eq 1.1-ext}) in $\R^N\setminus\Omega$,
 then
 $$\liminf_{|x|\to+\infty} u_0(x) |x|^{-\tau_-}>0.$$
 \end{proposition}
 {\bf Proof. } Let $f$ be a  nonnegative,  $C^1$ function  such that  $f\leq Qu_0^p$, $f$ has compact support.
By the comparison principle in the weak sense, we have that
$$u_0\geq u_f,$$
where $u_f$ is the classical solution of
  \begin{equation}\label{eq 2.1-green-ex}
\left\{\arraycolsep=1pt
\begin{array}{lll}
\cL_{\mu}^s u=  f  \quad\        {\rm in}\ \,   \R^N\setminus \bar \Omega , \\[2mm]
\quad \,  u= 0\quad  \     {\rm in}\ \     \bar\Omega,\\[2mm]
 \displaystyle\liminf_{|x|\to+\infty}u(x)|x|^{-\tau_+}=0.
  \end{array}
 \right.
\end{equation}
 By the strong maximum principle, we have that $u_f>0$ in $\R^N\setminus \bar \Omega$.

Now we can set that for $\varepsilon_0>0$ and $r_1>r_2>R_0$,
  $$f\geq \varepsilon_0\quad {\rm in}\ \ B_{r_1}\setminus B_{r_2},$$
 since $ Qu_0\geq  Qu_f^p>0$ a.e. in $\R^N\setminus \bar \Omega$.

Let
$$
  v^\sharp(x)=|x|^{2s-N}u_f(\frac{x}{|x|^2})\quad{\rm for}\  x\in O^\sharp,
 $$
then by (\ref{kt-1}) we have that
$$\cL_{\mu}^s v^\sharp(x)=|x|^{-2s-N}  f(\frac{x}{|x|^2})\quad{\rm in}\ \ B_{\frac1{R_0}}\setminus\{0\}.$$
Then it follows from Proposition \ref{lm 2.1} that for some $c_7>0$
$$v^\sharp(x)\geq c_7|x|^{\tau_+}\quad{\rm in}\ \ B_{\frac1{R_0}}\setminus\{0\},$$
which implies that
$$u_f(\frac{x}{|x|^2})\geq c_7|x|^{N-2s+\tau_+}\quad{\rm in}\ \ B_{\frac1{R_0}}\setminus\{0\}$$
and then
$$u_0(x)\geq c_7|x|^{ \tau_-}\quad{\rm in}\ \ \R^N\setminus B_{ R_0}, $$
where $\tau_++\tau_-=2s-N$.
  \hfill$\Box$\medskip

\begin{lemma}\label{lm 2.1-ex} Let $2s+\theta>0$,   $\tau_0<0$ and   $\{\tau_j\}_j$ be the sequence generated by (\ref{2.1-jj}).

$(i)$ For
\begin{equation}\label{2.00-ex}
 p\in\left[1,\ 1+\frac{2s+\theta}{-\tau_0}\right),
\end{equation}
  then $\{\tau_j\}_j$ is an increasing sequence of numbers and for any $\bar \tau>\tau_0$
there exists $j_0\in\N $  such that
\begin{equation}\label{2.3}
 \tau_{j_0}\ge \bar \tau\quad {\rm and}\quad \tau_{j_0-1}<\bar \tau.
\end{equation}

$(ii)$ For
\begin{equation}\label{2.00}
 p\in\left(0,\ 1\right),
\end{equation}
  then $\{\tau_j\}_j$ is an increasing sequence of numbers and for any $\bar \tau\in (\tau_0,\frac{2s+\theta}{1-p})$,
there exists $j_0\in\N $  such that
\begin{equation}\label{2.3}
 \tau_{j_0}\ge \bar \tau\quad {\rm and}\quad \tau_{j_0-1}<\bar \tau.
\end{equation}
\end{lemma}
 \noindent{\bf Proof. }
For $p\in(0,1+\frac{2s+\theta}{-\tau_0})$, we have that
$$\tau_1-\tau_0=2\alpha+\tau_0(p-1)>0$$
and
\begin{equation}\label{2.2}
\tau_j-\tau_{j-1} = p(\tau_{j-1}-\tau_{j-2})=p^{j-1} (\tau_1-\tau_0),
\end{equation}
which implies that the sequence $\{\tau_j\}_j$ is  increasing.
If $p\ge1 $, our conclusions are straightforward. If $p\in(0,1)$,
  in the case when $\tau_1\ge0$, we are done, and in the case when $\tau_1<0$,
it follows from (\ref{2.2}) that
\begin{eqnarray*}
\tau_j  &=&  \frac{1-p^j}{1-p}(\tau_1-\tau_0)+\tau_0\\
&\to&\frac{1}{1-p}(\tau_1-\tau_0)+\tau_0=\frac{2s+\theta}{1-p}\quad\  {\rm as}\ \, j\to+\infty,
\end{eqnarray*}
 then there exists $j_0>0$ satisfying (\ref{2.3}).
\hfill$\Box$\medskip

 \noindent {\bf Proof of Theorem \ref{teo 2}. }   By contradiction, we assume that $u_0$ is a positive super-weak solution of (\ref{eq 1.1-in}) in $\R^N \setminus\bar \Omega$.
From Proposition \ref{lm 2.1-ex}, we have that
$$u_0(x)\geq d_0 |x|^{\tau_-}\quad {\rm in}\ \, \R^N\setminus   \bar B_r$$
for $r>R_0$ and $d_0>0$.   It is worth noting that for $p\in(0, 1)$,
\begin{equation}\label{4.1-ex}
\frac{2s+\theta}{1-p}>\tau_-,
\end{equation}
thanks to $p< p^*_{\theta,\mu}$.

\smallskip

{\bf Part 5: $\mu> \mu_0$,  $\theta>\tau_--2s$ and $p\in(0, p^*_{\mu,\theta})$. } Let $\tau_0=\tau_-$,
 which verifies that for $x\in \R^N\setminus  \bar B_r$
$$
\cL^s_\mu  u_0(x)  \geq    q_\infty d_0^p |x|^{p\tau_-+\theta}=q_\infty d_0^p |x|^{\tau_1-2s}\quad {\rm in}\ \, \R^N\setminus  \bar B_r,
$$
where    $q_\infty>0$ and
$$\tau_1:=p\tau_0+\theta+2s.$$

If $p\tau_0+\theta\geq   \tau_+-2s$, then
$$ Qu_0^p\geq q_\infty d_0^p  |x|^{\tau_+-2s}$$
and a contradiction follows by Corollary \ref{re 2.1}. We are done.  \smallskip

If not, by Lemma \ref{teo 2.2}, we have that
$$u_0(x)\geq d_1|x|^{\tau_1}\quad {\rm in}\ \, \R^N\setminus   \bar B_r.$$
Iteratively, we recall that
$$\tau_j:=p\tau_{j-1}+\theta+2s,\quad j=1,2,\cdots .$$
Note that for $p\in(0,   p^*_{\mu,\theta})$
$$\tau_1-\tau_0=(p-1)\tau_0+\theta+2s >0.$$
 If $\tau_{j+1}=\tau_jp+\theta+2s\in  (\tau_-,\tau_+)$,
it following by Theorem \ref{teo 2.2}  that
$$u_0(x)\geq d_{j+1} |x|^{\tau_{j+1}},$$
where
$$\tau_{j+1}=p\tau_j+2s+\theta>\tau_j.$$
  If $p\tau_{j+1}+\theta\geq  \tau_-$, we are done by lemma \ref{teo 2.2}.  In fact, this iteration will stop after a finite number of iterations due to (\ref{4.1-ex}) for $p\in(0,1)$
  and $\tau_j\to+\infty$ as $j\to+\infty$ if $p\geq 1$.\smallskip

  {\bf Part  6: $\mu>\mu_0$,  $\theta>-2s$ and $p=p^*_{\mu,\theta}>1$. }
  Note  that for some $\sigma_0>0$
 $$\frac12 Qu_0^{p-1}(x)\geq \sigma_0|x|^{-2s} \quad{\rm in}\ \, \R^N\setminus B_r.$$
 Here we can assume that $\mu-\sigma_0\geq \mu_0$, otherwise, we take a smaller value for $\sigma_0$.
Therefore, we can write problem (\ref{eq 1.1}) as follows
\begin{equation}\label{eq 1.1 new-ex}
 \cL_{\mu-\sigma_0}^s u_0\geq \frac12 Q(x)   u_0^p \quad {\rm in}\ \  \R^N\setminus B_r,
 \end{equation}
 with the critical exponent
 $$p^*_{\mu-\sigma_0,\theta}=  1+\frac{2s+\theta}{-\tau_-(s,\mu-\sigma_0)} >    1+\frac{2s+\theta}{-\tau_+(s,\mu)}=p^*_{\mu,\theta}, $$
 where $\mu\in(\mu_0,+\infty)\mapsto \tau_+(s,\mu)$ is strictly increasing.
Thus a contradiction comes from {\bf part 5 } for (\ref{eq 1.1 new-ex}).\smallskip

  {\bf Part 7: $\mu>\mu_0$, $\theta=-2s$ and $p=1$. } This reduces the nonhomogeneous  problem
$$(-\Delta)^s u+(\mu+\mu')|x|^{-2s} u =f \quad{\rm in}\ \, B_r\setminus\{0\},$$
where $r>0$,  $f(x)\geq Q(x)-\mu'|x|^{-2s}\geq0$ and $\mu+\mu'<\mu_0$. While
 the above problem has no positive weak solutions for
 $\mu+\mu'<\mu_0$ from Corollary \ref{re 2.2}. \smallskip

{\bf Part  8:   $\mu=\mu_0$, $\theta>\tau_--2s$ and $p\in(0, p^*_{\mu,\theta}]$. } Since
  $$p\tau_0+\theta\geq  \tau_0,$$ then
$$ Qu_0^p\geq q_\infty d_0^p |x|^{\tau_+-2s}\quad{\rm for}\, \ |x|>r$$
and a contradiction follows by Corollary \ref{re 2.1}.   \hfill$\Box$

\setcounter{equation}{0}
\section{The Existence in the critical cases}

For $\mu\geq \mu_0$,
 let $\tau\in(-N,2\s)$,
$$v_\tau(x)=|x|^{\tau},\quad\forall\, x\in\R^N\setminus\{0\}$$
and a positive integer $m$,
$$w_{\tau,m}(x)=|x|^{\tau}(-\ln |x|)^m,\quad\forall\, x\in\R^N\setminus\{0\}. $$
%Note that for $r\in(0,\frac12)$,
%$$w_m'(r)=\tau_+ r^{\tau_+-1}(-\log r)^m-mr^{\tau_+-1}(-\log r)^{m-1} $$
%and
%$$w_m''(r)=\tau_+(\tau_+-1) r^{\tau_+-2}(-\log r)^m-m[2\tau_++N-2]r^{\tau_+-2}(-\log r)^{m-1}-m(m-1)r^{\tau_+-2}(-\log r)^{m-2},  $$
%then for $\mu>\mu_0$
% \begin{eqnarray*}
 % -\Delta w_m+\frac{\mu}{|x|^2} w_m&=& c_1|x|^{\tau_+-2}(-\log |x|)^{m-1}-c_2|x|^{\tau_+-2}(-\log |x|)^{m-2},
 %  \end{eqnarray*}
 %  where
  % $$c_1=m[2\tau_++N-2]>0 \quad{\rm and}\quad c_2=m(m-1).$$
It is known from  \cite[Lemma 2.1,\, 2.3]{CW} that
for $\tau \in (-N,2\s)$,  there holds that
$$
 (-\Delta)^\s v_\tau(x) = c_\s(\tau)|x|^{\tau-2\s}, \quad\ \forall\, x\in \R^N\setminus\{0\}$$
  with
\begin{equation}\label{2.1}c_\s(\tau) = 2^{2\s} \frac{\Gamma(\frac{N+\tau}{2})\Gamma(\frac{2\s-\tau}{2})}{\Gamma(-\frac{\tau}{2})\Gamma(\frac{N-2\s+\tau}{2})}.
\end{equation}
The function $c_\s: (-N, 2 \s) \mapsto \R $
is strictly concave and uniquely maximized at the point $\frac{2\s-N}{2}$ with the maximal value $2^{2\s}  \frac{\Gamma^2(\frac{N+2\s}4)}{\Gamma^2(\frac{N-2\s}{4})}.$

Moreover,
\begin{equation}
  \label{eq:c-symmetry}
c_\s(\tau)= c_\s(2\s-N-\tau) \qquad  {\rm for}\ \ \tau \in (-N,2\s)
\end{equation}
and
\begin{equation}
  \label{eq:c-asymptotics}
\lim_{\tau \to -N}c_\s(\tau) = \lim_{\tau \to 2 \s}c_\s(\tau)=-\infty.
\end{equation}

For $\mu>\mu_0$, $-N<\tau_-<\tau_+<2s$ are zero points of
$$c_\s(\tau)+\mu=0.$$

\begin{lemma}\label{lm}
For any positive integer $m$, we have that
 $$\cL_{\mu}^sw_{\tau_{\pm},m}(x)=|x|^{\tau_\pm-2s} \Big(\sum^m_{i=1} c_s^{(i)}(\tau_\pm) (\ln |x|)^{m-i}\Big),$$
where
$c_s^{(i)}(\tau_\pm)$ is the $i$-th-derivative of $c_s$ at $\tau=\tau_\pm$.

Furthermore,
$$c_s^{(i)}(\tau_+)<0\quad{\rm if} \ \, i\ {\rm is\ odd},\qquad c_s^{(i)}(\tau_+)>0\quad{\rm if} \ \, i\ {\rm is\ even}$$
and
$$c_s^{(i)}(\tau_-)>0\quad{\rm if} \ \, i\ {\rm is\ odd},\qquad c_s^{(i)}(\tau_-)<0\quad{\rm if} \ \, i\ {\rm is\ even}.$$
\end{lemma}
{\bf Proof.} Fix $ x\in \R^N\setminus\{0\}$ and let
$$F(x,\tau)=(-\Delta)^s v_\tau(x)+\frac{\mu}{|x|^{2s}} v_\tau(x)-c_\s(\tau)|x|^{\tau-2s}\equiv 0 \quad{\rm for}\ \, \tau\in (-N,2s),$$
then
\begin{eqnarray*}
 0&=& \left.\frac{\partial^{m} }{\partial \tau^m} F(x,\tau)\right|_{\tau=\tau_\pm}
 \\[2mm]&=& (-\Delta)^s w_{\tau_\pm,m}(x)+\frac{\mu}{|x|^{2s}}  w_{\tau_\pm,m}(x)-\sum^m_{i=0} c_s^{(i)}(\tau_\pm)|x|^{\tau_{\pm}-2s}( \ln |x|)^{m-i}
  \\[2mm]&=&\cL_{s,\mu}w_{\tau_\pm,m}(x)-\sum^m_{i=1} c_s^{(i)}(\tau_\pm)|x|^{\tau_\pm-2s}(\ln |x|)^{m-i},
   \end{eqnarray*}
where we use the fact that $c_s^{(0)}(\tau_+)=c_s(\tau_+)=0$.

Note that
\begin{eqnarray*}
c_\s(\tau) |x|^{\tau-2\s}  &=&   (-\Delta)^\s v_\tau(x)
\\[2mm] &=& -\frac{C_{N,\s}}2 \int_{\R^N}\frac{|x+y|^{\tau}+|x-y|^\tau-2|x|^\tau}{|y|^{N+2\s}}\,dy\\[2mm]
    &=& -\frac{C_{N,\s}}2|x|^{\tau-2\s}\int_{\R^N}\frac{|e_1+z|^{\tau}+|e_1-z|^\tau-2}{|z|^{N+2\s}}\,dz \quad  {\rm for}\  x \in \R^N \setminus \{0\},
\end{eqnarray*}
where $e_1=(1,0,\cdots,0)\in\R^N$. Thus
$$c_\s (\tau)=-\frac{C_{N,\s}}2\int_{\R^N}\frac{|x-e_1|^\tau+|x+e_1|^\tau-2}{|x|^{N+2\s}}\,dx.$$

Consequently, for $\tau \in (-N,2\s)$, we have for positive integers $n$ that
\begin{equation}\label{cs2n}
c_\s^{(2n)}(\tau)=-\frac{C_{N,\s}}{2}\int_{\R^N}\frac{|e_1-x|^{\tau} (\log|e_1-x|)^{2n}+|e_1+x|^{\tau}(\log|e_1+x|)^{2n}}{|x|^{N+2\s}}dx<0,
\end{equation}
then $c_\s^{(2n-2)}(\tau_+)$ is strictly concave,  and by the symmetric property (\ref{eq:c-symmetry}), we have that
$c_\s^{(2n-1)}(\cdot)$ is strictly decreasing at $\tau=\tau_+$ and strictly increasing at  $\tau=\tau_-$, that is
$$c_\s^{(2n-1)}(\tau_+)<0<c_\s^{(2n-1)}(\tau_-).$$
This completes the proof. \hfill$\Box$

\subsection{ The existence in the sublinear critical case    }
\noindent{\bf Proof of Theorem \ref{teo cri ex 1}. }  For $\tau_+-2s<\theta< -2s$, we see that $p^\#_{\mu,\theta}\in(0,1)$.

Take an even integer $m\geq \frac1{1-p^\#_{\mu,\theta}}$,  then $w_{\tau_+,m}$ is nonnegative and there exists $r_m\in(0,\frac1e)$ such that for $0<|x|<r_m$,
\begin{eqnarray*}
\cL_{\mu}^sw_{\tau_+,m}(x)&=&|x|^{\tau_+-2s} \Big(\sum^m_{i=1}c_s^{(i)}(\tau_+) (-\ln |x|)^{m-i}\Big)\\&\geq& \frac12(-c_s'(\tau_+)) |x|^{\tau_+-2s} (-\ln |x|)^{m-1},
\end{eqnarray*}
where $-c_s'(\tau_+)>0$ and $m-1$ is odd. On the other hand, by assumption (\ref{pt r1-l}) and the choice of $m$
$$Q(x)w_{\tau_+,m}(x)^{p^\#_{\mu,\theta}} \leq q_0 |x|^{\tau_+-2s} (-\ln |x|)^{m p^\#_{\mu,\theta}}\leq q_0 |x|^{\tau_+-2s} (-\ln |x|)^{m-1}.$$

For the existence with $\Omega=B_{r_m}$,  for any  $t\in\Big(0,\, (\frac{2q_0}{-c_s'(\tau_+)})^{\frac1{1-p^\#_{\mu,\theta}}}\Big]$,  function
$tw_{\tau_+,m}$  verifies (\ref{eq 1.1-in}) in $B_{r_m}\setminus\{0\}$.\hfill$\Box$\bigskip

 \noindent{\bf Proof of Theorem \ref{teo cri ex 2}. }  For $\tau_--2s<\theta< -2s$, computations show that $p^*_{\mu,\theta}\in(0,1)$.

Take an even integer $m\geq \frac1{1-p^*_{\mu,\theta}}$,  then $w_{\tau_-,m}$ is nonnegative and there exists $R_m>e$ such that for $|x|>R_m$,
$$\cL_{\mu}^s w_{\tau_-,m}(x)=|x|^{\tau_--2s} \Big(\sum^m_{i=1}c_s^{(i)}(\tau_-) (\ln |x|)^{m-i}\Big)\geq \frac12c_s'(\tau_-) |x|^{\tau_+-2s} (\ln |x|)^{m-1},$$
where $c_s'(\tau_-)>0$. By assumption (\ref{pt r1-l}) and the choice of $m$
$$Q(x)w_{\tau_-,m}(x)^{p^*_{\mu,\theta}} \leq q_0 |x|^{\tau_--2s} (\ln |x|)^{m p^*_{\mu,\theta}}\leq q_0 |x|^{\tau_--2s} (\ln |x|)^{m-1}.$$

For the existence,   setting $\Omega=\R^N\setminus \bar B_{R_m}$,   for any  $t\in\Big(0,\, (\frac{2q_0}{c_s'(\tau_-)})^{\frac1{1-p^*_{\mu,\theta}}}\Big]$,  function
$tw_{\tau_-,m}$  verifies (\ref{eq 1.1-ext}) in $\R^N\setminus \bar B_{R_m}$.\hfill$\Box$

 \bigskip\bigskip

   \noindent{\bf \small Acknowledgements:} {\footnotesize This work is is supported by NNSF of China, No: 12071189 and 12001252,
by the Jiangxi Provincial Natural Science Foundation, No: 20202BAB201005,
by the Science and Technology Research Project of Jiangxi Provincial Department of Education,
No: 200307 and 200325.}

\end{document}